\documentclass [12pt]{article}
\usepackage {amsfonts}
\usepackage {mathrsfs}
\usepackage {amsmath}
\usepackage {latexsym}
\usepackage {amssymb}
\usepackage {amsthm}
\usepackage {amscd}
\date{}
\title{\bf Basis of Nichols Braided Lie Algebras}
\author{\small  Weicai Wu $^{a}$,    Jing Wang $^{b}$ and Shouchuan Zhang $^{c}$  \\
\small $a$. Department  of Mathematics,    Zhejiang  University\\
\small   Hangzhou  310007,    \ P.R. China \\
\small $b$. Department of mathematics,  School of Science,Beijing Forestry University, \\
\small    Beijing, 100083  \ P.R. China \\
\small $c$. Department  of Mathematics,    Hunan University\\
\small   Changsha  410082,    \ P.R. China \\
\small {\tt Emails: weicaiwu@hnu.edu.cn(WC); wang\_jing619@163.com (JW)}\\
 }

\begin{document}
\newtheorem{Proposition}{Proposition}[section]
\newtheorem{Theorem}[Proposition]{Theorem}
\newtheorem{Definition}[Proposition]{Definition}
\newtheorem{Corollary}[Proposition]{Corollary}
\newtheorem{Lemma}[Proposition]{Lemma}
\newtheorem{Example}[Proposition]{Example}
\newtheorem{Remark}[Proposition]{Remark}

\maketitle %\addtocounter{section}{-1}

\begin {abstract} Assume that $V$ is  a  braided vector space with      diagonal type. It is shown that  a monomial belongs to Nichols braided Lie algebra $\mathfrak L(V)$ if and only if
this monomial is connected. A basis of Nichols braided Lie algebra and dimension of Nichols braided Lie algebra of finite Cartan type
 are obtained.

\end {abstract}

\section {Introduction}\label {s0}
The question of finite-dimensionality of Nichols algebras dominates an important part of the recent developments in the theory of (pointed) Hopf algebras(see e.g. \cite {AHS08,  AS10,  He05,  He06a,  He06b,  WZZ15a,  WZZ15b}. The interest in this problem comes from the lifting method by Andruskiewitsch and Schneider to classify finite dimensional (Gelfand-Kirillov) pointed Hopf algebras,   which are generalizations of quantized enveloping algebras of semi-simple Lie algebras.

The classification of arithmetic root systems is obtained in \cite {He05} and \cite {He06a}.
It is proved in \cite {WZZ15b} that   Nichols algebra $\mathfrak B(V)$  is
finite-dimensional if and only if Nichols braided Lie algebra $\mathfrak L(V)$ is finite-dimensional.

The theory of Lie superalgebras has been developed systematically,
which includes  the representation theory and classifications of
simple Lie superalgebras and their varieties \cite {Ka77} \cite
{BMZP92}.  In many physical applications or in pure mathematical
interest,  one has to consider not only ${\bf Z}_2$- or ${\bf Z}$-
grading but also $G$-grading of Lie algebras,  where $G$ is an
abelian group equipped with a skew symmetric bilinear form given
by a 2-cocycle. Lie algebras in symmetric and more general
categories were discussed in \cite {Gu86} and \cite {GRR95}.
A sophisticated multilinear version of the Lie bracket was
considered in \cite {Kh99} \cite {Pa98}. Various generalized Lie
algebras have already appeared  under different names,  e.g. Lie
color algebras,  $\epsilon $ Lie algebras \cite {Sc79},  quantum and
braided Lie algebras,   generalized Lie algebras \cite {BFM96},
$H$-Lie algebras \cite {BFM01} and braided $m$-Lie algebras \cite {ZZ03}.

In this paper it is shown that  a monomial belongs to Nichols braided Lie algebra $\mathfrak L(V)$ if and only if
this monomial is connected. A basis of Nichols braided Lie algebra and dimension of Nichols braided Lie algebra of finite Cartan type
 are obtained for a  braided vector space $V$ with      diagonal type.

This paper is organized as follows. In Section 1 we provide some preliminaries and set our notations. In  Section 2 we obtain that a monomial belongs to $\mathfrak L(V)$ if and only if
this monomial is connected when $V$ is  a  braided vector space with      diagonal type. In  Section 3 we obtain a basis of Nichols braided Lie algebra.
In Section 4 we obtain  dimension of Nichols braided Lie algebra $\mathfrak L(V) $
 with finite Cartan type. In  Section 5 we find some non-zero monomials in  Nichols algebras. In  Section 6 we give relations between connected components of    graphes  and Nichols Braided Lie Algebras.

\section*{ Preliminaries}\label {s1}

Let $\mathfrak B(V)$ be the Nichols algebra generated by vector space $V$. Throughout this paper  braided vector space $V$ is of diagonal type with basis
$x_1,  x_2,  \cdots,  x_n$ and $C(x_i \otimes x_j) = p_{ij} x_j \otimes x_i$,  $n>1$ without special announcement.
Define  linear map $p$ from $\mathfrak B(V) \otimes \mathfrak B(V) $ to $F$ such that
$p(u\otimes v) = \chi (deg (u),   deg (v)),  $ for any homogeneous element $u,   v \in \mathfrak B(V).$
For convenience,   $p(u\otimes v)$ is denoted by $p_{uv}$.  Let $\widetilde{p}_{uv} := p_{uv}p_{vu}.$

Denote ${\rm ord } (p_{uu})$ the  order of $p_{uu}$ with respect to multiplication. Let $|u|$ denote length of  homogeneous element $u\in \mathfrak B(V).$
Let $D =: \{[u] \mid [u] \hbox { is a hard super-letter }\}$, $\Delta ^+(\mathfrak B(V)): =  \{ \deg (u) \mid [u]\in D\}$,
$ \Delta (\mathfrak B(V)) := \Delta ^+(\mathfrak B(V)) \cup \Delta ^-(\mathfrak B(V))$,
which is called the root system of $V.$
 If $ \Delta (\mathfrak B(V))$ is finite,
then it is called an arithmetic root system. Let  $\mathfrak L(V)$ denote the  braided Lie algebras generated by $V$ in $\mathfrak B(V)$ under Lie operations $[x, y]=yx-p_{yx}xy$,  for any homogeneous elements $x, y\in \mathfrak B(V)$. $(\mathfrak L(V), [\ ])$ is called Nichols braided Lie algebra of $V$. The other notations are the same as in \cite {WZZ15a}.

\vskip.1in
Recall the dual $\mathfrak B(V^*) $ of Nichols algebra $\mathfrak B(V) $ of rank $n$
in \cite [Section 1.3]{He05}. Let $y_{i}$ be a dual basis of $x_{i}$. $\delta (y_i)=g_i^{-1} \otimes y_i$,  $g_i\cdot y_j=p_{ij}^{-1}y_j$ and $\Delta (y_i)=g_i^{-1}\otimes y_i+y_i\otimes 1.$ There exists a bilinear map
$<\cdot, \cdot >$: $(\mathfrak B(V^*)\# FG)\times\mathfrak B(V)$ $\longrightarrow$ $\mathfrak B(V)$ such that

$<y_i, uv>=<y_i, u>v +g_i^{-1}.u<y_i, v>$ and $<y_i, <y_j, u>>=<y_iy_j, u>$

\noindent for any $u, v\in\mathfrak B(V)$. Furthermore,  for any $u\in \oplus _{i=1}^\infty \mathfrak B(V)_{(i)}$,
one has that  $u=0$ if and only if $<y_i, u>=0$ for any $1\leq i \leq n.$

We have the braided Jacobi identity as follows:
\begin {eqnarray}\label {e0.1} [[u, v], w]=[u, [v, w]]+p_{vw}^{ -1} [[u, w], v]+(p_{wv}-p_{vw}^{-1})v\cdot[u, w].
\end {eqnarray}
\begin {eqnarray}\label {e0.2} [u, v \cdot w]=p_{wu}[uv]\cdot w + v \cdot [uw].\end {eqnarray}

$\Gamma (V)$ be the  generalized Dynkin diagram of $V$ omitted $p_{x_i,  x_i},\widetilde{p}_{x_i,  x_j}$ for $1\le i,   j \le n $,  which is called pure  generalized Dynkin diagram of $V$.

Let   $u=  h_1 h_2\cdots  h_m$ be   a  monomial with $ h_j = x_{i_j}$  for $1\leq i_1,i_2,\ldots,i_m\leq n,1\le j \le m$. ${\rm deg} (u) = \lambda _1 e_1 + \cdots + \lambda _n e_n, $ where $ {\rm deg }(x_i) = e_i.$ Let ${\rm deg}_{x_i} (u) := \lambda _i$. Let $\mu (u):= \{ x_{i_1},  \cdots,  x_{i_m}\}$  and  $\Gamma (u)$ be a pure generalized Dynkin sub-diagram generated by $\mu (u)$.
 If $\Gamma (u)$
is connected and  $p_{h_i,   h_i} \not= 1$ for $1\le i \le m$,  then $u$ is called connected (or $\mu (u) $ is called connected).
If $ h_i  h_{i+1} = p_{ h_{i+1} h_i}  h_{i+1} h_i $,  then $ h_1  h_2\cdots  h_m = p_{ h_{i+1} h_i} h_1 h_2 \cdots  h_{i+1}  h_i \cdots  h_m$  and $  h_1 h_2 \cdots  h_{i+1}  h_i \cdots  h_m$ is called an elementary quantum equivalent of  $ h_1 h_2 \cdots  h_{i+1}  h_i \cdots  h_m $. If $v_{1}$ is a monomial and $v _{j+1}$ is  an elementary quantum equivalent of $v_j$   for $1\le j \le r-1$,  then $v_r$ is called  a quantum equivalent of $v_1, $ written  as $v_1 \sim v_r.$ This is an equivalent relation since $ h_i  h_{i+1} = p_{ h_{i+1} h_i}  h_{i+1} h_i $ if and only if $p_{ h_{i+1} h_i}p_{ h_{i} h_{i+1}} =1.$
If there exist $x_i \in \mu (u)$ and   $x_j \in \mu (v)$   such that $\widetilde{p} _{x_i,  x_j} \not=1$,  then we say that it is connected between monomial $u$ and monomial $v$, written $u\circ v$ in short.

Throughout,  $\mathbb Z =: \{x \mid  x \hbox { is an integer}\}.$ $\mathbb R =: \{ x \mid x \hbox { is a real number}\}$.
$\mathbb N_0 =: \{x \mid  x \in \mathbb Z, x\ge 0\}.$
$\mathbb N =: \{x \mid  x \in \mathbb Z,  x>0\}$.  $F$ denotes the base field,   which is an algebraic closed field with  characteristic zero. $F^{*}=F\backslash\{0\}$. $S_{n}$ denotes symmetric group, $n\in\mathbb N$. For any set $X$, $\mid X\mid$ is the cardinal of $X$. $int(a)$ means the biggest integer not greater than $a\in \mathbb R$.

\section{The structure of   Nichols braided Lie algebras}\label {s1}
In this section we obtain that a monomial belongs to $\mathfrak L(V)$ if and only if
this monomial is connected when $V$ is  a  braided vector space with      diagonal type.

\begin{Lemma}\label{2.1} Assumed that $u, v, w$ are homogeneous elements in $\mathfrak L(V)$. Let $a,  b,  c,  d,  e,  f$ denote  $1 - p_{wv}p_{vw},  1 -p_{uw}p_{wu},  1 - p_{uv}p_{vu},  1 - p_{uv}p_{vu}p_{uw}p_{wu},  1 - p_{uv}p_{vu}p_{wv}p_{vw},  1 - p_{wv}p_{vw}p_{uw}p_{wu}$, respectively. If $\mid\{r|r\in\{uv,uw,vw\} \hbox { and }r\in\mathfrak L(V)\}\mid \geq 2$ and  $ \mid \{t|t\in\{a, b, c\} \hbox { and } t\neq0\} \mid \ge  1$,  then
 $uvw, uwv, vwu, vuw, wuv, wvu\in\mathfrak{L}(V)$.
\end {Lemma}
\noindent {\it Proof.} Without loss of generality, we let $uv,uw\in\mathfrak L(V)$. If $vuw\notin\mathfrak L(V)$, then $e=0,f=0$ by
\cite [Lemma 3.2(i)] {WZZ15b} and \cite [Lemma 4.12] {WZZ15a}. If $a=0$, then $b=0,c=0,$  which  is a contradiction. If $a\neq0$, then $b\neq0,c\neq0$, one obtains a contradiction to \cite [Lemma 3.1] {WZZ15b}.  Consequently,   $vuw\in\mathfrak L(V)$. Similarly,  we can obtain others.  \hfill $\Box$

\begin {Lemma} \label {2.3} Assume $h_i \in \{x_1,  \cdots,  x_n\}$ for $1\le i \le m.$

{\rm (i)} If it is disconnected between monomial $u$ and monomial $v$ (i.e. $\widetilde{p} _{x_i,  x_j} =1$  for any  $x_i \in \mu (u)$,  $x_j \in \mu (v)$ ),  then $[u,  v] =0.$

{\rm (ii)} If $\mu(h_{1}h_{2}\cdots h_{m})$ is  disconnected,  then $\sigma (h_{1},  h_{2},  \cdots,   h_{m})=0$ for any  method $\sigma$ of adding bracket on $h_1,  h_2,  \cdots,  h_m.$.

{\rm (iii)} If $h_{1}h_{2}\cdots h_{m}\neq0$ and $\mu(h_{1}h_{2}\cdots h_{m})$ is disconnected,  then $h_{1}h_{2}\cdots h_{m}\notin\mathfrak{L}(V)$.

\end {Lemma}

\noindent {\it Proof.} {\rm (i)}  $u$ and $v$ are quantum commutative
(i.e. $uv= p_{u, v}vu$ ) since $x_i x_j = p_{x_i,  x_j}x_jx_i$ for any  $x_i \in \mu (u)$,  $x_j \in \mu (v)$.

{\rm (ii)} We show this by induction on $m$. $[h_2,  h_1]= [h_1,  h_2] =0$ for $m=2.$  For $m>2,  $ $\sigma (h_{1},  h_{2},  \cdots,  h_{m}) = [\sigma _1 (h_{\tau_{(1)}}h_{\tau_{(2)}}\cdots h_{\tau_{(t)}}),  \sigma _2 ( h_{\tau_{(t+1)}}h_{\tau_{(t+2)}}\cdots h_{\tau_{(m)}})],\tau\in \mathbb S_{m}$. If both $h_{\tau_{(1)}}h_{\tau_{(2)}}\cdots h_{\tau_{(t)}}$ and $h_{\tau_{(t+1)}}h_{\tau_{(t+2)}}\cdots h_{\tau_{(m)}}$ are connected,  then it is disconnected between  $h_{\tau_{(1)}}h_{\tau_{(2)}}\cdots h_{\tau_{(t)}}$ and $h_{\tau_{(t+1)}}h_{\tau_{(t+2)}}\cdots h_{\tau_{(m)}}$. By Part {\rm (i)},   $\sigma (h_{1},  h_{2},  \cdots,  h_{m})=0$.
If either  $h_{\tau_{(1)}}h_{\tau_{(2)}}\cdots h_{\tau_{(t)}}$ or $h_{\tau_{(t+1)}}h_{\tau_{(t+2)}}\cdots h_{\tau_{(m)}}$ is  disconnected,
then either  $\sigma _1 (h_{\tau_{(1)}}h_{\tau_{(2)}}\cdots h_{\tau_{(t)}})=0$ or $ \sigma _2 (h_{\tau_{(t+1)}}h_{\tau_{(t+2)}}\cdots h_{\tau_{(m)}})=0 $ by inductive hypothesis.

{\rm (iii)}  It follows from Part {\rm (ii)}.
\hfill $\Box$

\begin {Lemma} \label {p2.3}

(i) If $u$ is  a monomial, then there exist monomials $u_{1},u_{2},\ldots,u_{r}$ such that $u \sim  u_1u_2 \cdots u_r $ for $\{ \Gamma (u_1),  \Gamma ( u_2),  \cdots,  \Gamma(u_r) \}$ is complete set of a connected component of $\Gamma (u)$ with ${\rm deg}( u) =\sum \limits _{i=1} ^r {\rm deg}( u_i)$ and $\mid u_1 \mid\geq\mid u_2 \mid\geq\cdots\geq\mid u_r \mid$,  (which is called a  decomposition of  connected components  of $u$)

(ii)  If a monomial $u$ is connected with $\mid u \mid >1$,  then there exist two  connected  monomials $v$ and $w$ such that $v\circ w$  with $u \sim vw$.

\end {Lemma}
{\bf Proof.} (i)  If $u$ is connected, it is clear. If $u$ is not connected, we show it by induction on $\mid u \mid $. It is clear when  $\mid u \mid =1.$  Now assume $\mid u \mid >1$. Let $\Omega := \{ v \mid
    v \hbox { is connected   }  \}$.  Let $ u = v_1v_2v_3$ and $v_2$ be in $\Omega$ such that $\mid v_2 \mid = max \{ \mid v \mid \ \  \mid  v \in  \Omega \}$.

If $v_3\in F^*$, then $v_2v_3 \sim v_3v_2$. If $v_2v_3 \neq kv_3v_2$ for $\forall\ k\in F^*$, then there exists $t$ such that $ v_2\circ h_{i_t}$ with
$v_3 = h _{i_1} \cdots  h_{i_r} $  and  $v_2  h_{i_s} \sim  h_{i_s } v_2$ for $1\le s < t$. Consequently,  $u \sim w_1 v_2 h_{i_t}w_3 $ and $v_2 h_{i_t}$ is connected. which is a contradiction. Consequently, $v_2v_3 \sim v_3v_2$.

If $v_1\in F^*$, then $v_1v_2 \sim v_2v_1$. If $v_1v_2 \neq kv_2v_1$ for $\forall\ k\in F^*$, then there exists $t$ such that $ h_{i_t}\circ v_2$ with
$v_1 = h _{i_1} \cdots  h_{i_r} $  and  $v_2  h_{i_s} \sim  h_{i_s } v_2$ for $t+1< s \leq r$. Consequently,  $u \sim w_1 h_{i_t}v_2 w_3 $ and $h_{i_t}v_2 $ is connected. which is a contradiction. Consequently, $v_1v_2 \sim v_2v_1$.

Thus  $u =v_1 v_2 v_3 \sim v_2 v _1v_3 $ and $v_1v_3 \sim u_2u_3 \cdots u_r$ is a  decomposition of  connected components  of $v_1v_3$ since $\mid v_1v_3 \mid < \mid u \mid. $ Set $u_1:= v_2 $.

(ii) Let $ u = h_1 h_2\cdots  h_m$   and $ h_2\cdots  h_m \sim u_1u_2\cdots u_r$ be a decomposition of connected components of $ h_2\cdots  h_m$.
Consequently,  $ h_1u_1\cdots u_{r-1}$ is connected since $ h_1\circ u_i$ for $1\le i \le r$.
\hfill  $\Box$

\begin {Theorem} \label {pp2.3}Assume $h_i \in \{x_1,  \cdots,  x_n\}$ and  $p_{ h_1, h_1 } \not= 1$ when $h_i = h_1$  for $1\le i \le m$.

(i)
 If $ h_{1}\cdots  h_{b},  h_{b+1}\cdots  h_{m}\in\mathfrak{L}(V)$ and
$ h_{1}\cdots  h_{b}\circ h_{b+1}\cdots  h_{m}$ with $h_{1}\cdots  h_{b}\not= 0$ or $ h_{b+1}\cdots  h_{m} \not= 0$,
then there exist $\tau\in \mathbb  S_m $ such that $ h_{1}\cdots  h_{m}\sim h_{\tau(1)}\cdots  h_{\tau(m)}$ with

\noindent $ h_{\tau(1)}\cdots h_{\tau(m-1)}\in\mathfrak{L}(V)$ and  $ h_{\tau(1)}\cdots h_{\tau(m-1)}\circ h_{\tau(m)}$,
or with $ h_{\tau(2)}\cdots  h_{\tau(m)}\in\mathfrak{L}(V)$ and $ h_{\tau(2)}\cdots  h_{\tau(m)}\circ h_{\tau(1)}$.

(ii)
If $  h_{1} h_{2}\cdots  h_{m-1}\in\mathfrak{L}(V) $  and $  h_{1} h_{2}\cdots  h_{m-1}\circ h_m$ or  $   h_{2} h_{3}\cdots h_{m}\in\mathfrak{L}(V) $  and $  h_{2} h_{3}\cdots  h_{m}\circ h_1, $ then $ h_{1} h_{2}\cdots  h_{m}\in\mathfrak{L}(V)$.

(iii) If $0\not=   h_{1} h_{2}\cdots  h_{m}\in\mathfrak{L}(V) $,  then there exists $\tau \in \mathbb S_m$ such that $ h_{1}\cdots  h_{m}\sim h_{\tau(1)}\cdots  h_{\tau(m)}$ with $ 0\neq h_{\tau (1)} h_{\tau(2)}\cdots h_{\tau(m-1)}\in\mathfrak{L}(V)$ and $ h_{\tau (1)} h_{\tau(2)}\cdots  h_{\tau(m-1)}\circ h_{\tau(m)}$,   or $0\neq h_{\tau (2)} h_{\tau(3)}\cdots  h_{\tau(m)}\in\mathfrak{L}(V)$ and $ h_{\tau (2)} h_{\tau(3)}\cdots  h_{\tau(m)}\circ h_{\tau(1)}$ .

 (iv)
If monomial $u$ is connected,  then $u \in  \mathfrak{L}(V). $

\end {Theorem}
{\bf Proof.}
We show (i),  (ii),  (iii) and (iv) by induction on the length of $| h_{1} h_{2}\cdots  h_{m}|=m$.

  Assume $m=2$. (i) and (iv) are clear. If $ h_{1}\neq  h_{2}$,  then (ii), (iii) follows from \cite [Lemma 4.12]{WZZ15a}, \cite [Lemma 5.2]{WZZ15a},  respectively. If $ h_{1}=  h_{2}$,  then (ii) and (iii) follow from \cite [Lemma 4.3]{WZZ15a}  and  \cite [Lemma 1.3.3(i)]{He05}.

Now  $m>2$.

(i) If $b=m-1$ or $b=1$,  let $\tau=id$. Now assume that $1<b<m-1$. There exist $1 \le a \le b$ and $b+1 \le c \le m$ such that $\widetilde{p} _{h_a,  h_c} \not= 1$.

(1). Assumed that $ h_{1}\cdots  h_{b}\neq0$,  we proceed by induction over $b$.

\noindent We know there exist $\tau\in \mathbb  S_{1, \ldots, b}$ such that $ h_{1}\cdots  h_{b}\sim  h_{\tau(1)}\cdots  h_{\tau(b)}$ with

\noindent (a) $0\neq h_{\tau(1)}\cdots h_{\tau(b-1)}\in\mathfrak{L}(V)$, $ h_{\tau (1)} h_{\tau(2)}\cdots  h_{\tau(b-1)}\circ h_{\tau(b)}$
or with

\noindent (b) $0\neq  h_{\tau(2)}\cdots h_{\tau(b)}\in\mathfrak{L}(V)$,  $ h_{\tau (2)} h_{\tau(3)}\cdots  h_{\tau(b)}\circ h_{\tau(1)}$ by $ h_{\tau(1)}\cdots  h_{\tau(b)}\in\mathfrak L( V)$ and induction hypotheses of (iii).

(a)$_{1}$.  If $ h_{b+1} \cdots  h_{m}\circ h_{\tau(b)}$,  then
$ h_{\tau(b)} h_{b+1}\cdots  h_{m}\in\mathfrak L( V)$ by induction hypotheses of (ii),  and $h_{\tau(1)}\cdots h_{\tau(b-1)}\circ h_{\tau(b)} h_{b+1}\cdots  h_{m}$ since
$ h_{\tau (1)} h_{\tau(2)}\cdots  h_{\tau(b-1)}\circ h_{\tau(b)}$. It is proved since
$| h_{\tau(1)}\cdots  h_{\tau(b-1)}|=b-1<b$ and induction hypotheses.

(a)$_{2}$.  If all $j_{1}\in\{b+1, \ldots, m\}$ such that $\tilde{p}_{ h_{\tau(b)},  h_{j_{1}}}=1$,  then
$ h_{\tau(b)} h_{b+1}\cdots  h_{m}\sim  h_{b+1}\cdots  h_{m} h_{\tau(b)}$. On the other hand,  we know $h_{\tau(1)}\cdots h_{\tau(b-1)}\circ  h_{b+1}\cdots  h_{m}$ since
$ h_{1}\cdots  h_{b}\circ h_{b+1}\cdots  h_{m}$,
then $ h_{\tau_{1}\tau(1)}\cdots  h_{\tau_{1}\tau(b-1)} h_{\tau_{1}(b+1)}\cdots h_{\tau_{1}(m)}\in\mathfrak L( V)$,  $ h_{\tau(1)}\cdots  h_{\tau(b-1)} h_{b+1}\cdots  h_{m}\sim  h_{\tau_{1}\tau(1)}\cdots h_{\tau_{1}\tau(b-1)} h_{\tau_{1}(b+1)}\cdots  h_{\tau_{1}(m)}$ for some $\tau_{1}\in \mathbb  S_{\tau(1), \ldots, \tau(b-1), b+1, \ldots, m} $ by the induction hypotheses of (i) and (ii). We know $h_{\tau_{1}\tau(1)}\cdots  h_{\tau_{1}\tau(b-1)} h_{\tau_{1}(b+1)}\cdots h_{\tau_{1}(m)}\circ h_{\tau(b)}$  since

\noindent $ h_{\tau (1)} h_{\tau(2)}\cdots  h_{\tau(b-1)}\circ h_{\tau(b)}$. We obtain
$ h_{1}\cdots  h_{m}
\sim  h_{\tau(1)}\cdots  h_{\tau(b)} h_{b+1}\cdots  h_{m}
$

\noindent $\sim  h_{\tau(1)}\cdots h_{\tau(b-1)} h_{b+1}\cdots  h_{m} h_{\tau(b)}\sim  h_{\tau_{1}\tau(1)}\cdots  h_{\tau_{1}\tau(b-1)} h_{\tau_{1}(b+1)}\cdots h_{\tau_{1}(m)} h_{\tau(b)}$

\noindent $=  h_{\tau_{1}\tau(1)} \cdots h_{\tau_{1}\tau(b-1)}$ $h_{\tau_{1}\tau(b+1)}\cdots  h_{\tau_{1}\tau(m)} h_{\tau_{1}\tau(b)}$ by $ h_{\tau_{1}(b+1)}\cdots  h_{\tau_{1}(m)}$
$= h_{\tau_{1}\tau(b+1)}\cdots  h_{\tau_{1}\tau(m)}$ since $\tau\in \mathbb  S_{1, 2, \ldots, b}$,
 and
$ h_{\tau_{1}\tau(b)}= h_{\tau(b)}$ since $\tau_{1}\in \mathbb  S_{\tau(1), \tau(2), \ldots, \tau(b-1), b+1, \ldots, m}$. Set
 $$\tau':={\tiny \left(       \begin{array}{ccccccc}
              1&\cdots&b-1&b&\cdots&m-1&m \\
               \tau_{1}\tau(1)&\cdots&\tau_{1}\tau(b-1)&\tau_{1}\tau(b+1)&\cdots& \tau_{1}\tau(m)&\tau_{1}\tau(b)
                   \end{array}\right) },  $$
                   then we obtain $ h_{\tau'(1)}\cdots h_{\tau'(m-1)}\in\mathfrak{L}(V), h_{\tau'(1)}\cdots h_{\tau'(m-1)}\circ h_{\tau'(m)}$ and

\noindent $ h_{1}\cdots  h_{m}\sim h_{\tau'(1)}\cdots  h_{\tau'(m)}$ for some $\tau'\in \mathbb  S_{1, \ldots, m}.$

(b)$_{1}$.  If $ h_{\tau(2)}\cdots h_{\tau(b)}\circ h_{b+1}\cdots  h_{m}$,  then  $ h_{\tau_{1}\tau(2)}\cdots  h_{\tau_{1}\tau(b)} h_{\tau_{1}(b+1)}\cdots h_{\tau_{1}(m)}\in\mathfrak L( V)$  and $ h_{\tau(2)}\cdots h_{\tau(b)} h_{b+1}\cdots  h_{m}
\sim  h_{\tau_{1}\tau(2)}\cdots
 h_{\tau_{1}\tau(b)} $ \ $h_{\tau_{1}(b+1)}\cdots  h_{\tau_{1}(m)}$ for some $\tau_{1}\in \mathbb  S_{\tau(2), \ldots, \tau(b), b+1, \ldots, m}$ by the induction hypotheses of (i) and (ii). We know $h_{\tau_{1}\tau(2)}\cdots  h_{\tau_{1}\tau(b)} h_{\tau_{1}(b+1)}\cdots h_{\tau_{1}(m)}\circ h_{\tau(1)}$ since $ h_{\tau (2)} h_{\tau(3)}\cdots  h_{\tau(b)}\circ h_{\tau(1)}$.
We obtain
$ h_{1}\cdots  h_{m}\sim  h_{\tau(1)}\cdots  h_{\tau(b)} h_{b+1}\cdots h_{m}\sim  h_{\tau(1)} h_{\tau_{1}\tau(2)}\cdots h_{\tau_{1}\tau(b)} h_{\tau_{1}(b+1)} $ $\cdots  h_{\tau_{1}(m)}= h_{\tau_{1}\tau(1)} h_{\tau_{1}\tau(2)}\cdots  h_{\tau_{1}\tau(b)} $ $h_{\tau_{1}\tau(b+1)}\cdots h_{\tau_{1}\tau(m)}$ by

\noindent $ h_{\tau_{1}(b+1)}\cdots h_{\tau_{1}(m)}= h_{\tau_{1}\tau(b+1)}\cdots  h_{\tau_{1}\tau(m)}
$ since
$\tau\in \mathbb  S_{ 1, \ldots, b}$ and $ h_{\tau_{1}\tau(1)}= h_{\tau(1)}$ since
$\tau_{1}\in \mathbb  S_{\tau(2), \ldots, \tau(b), b+1, \ldots, m}$. Set $\tau'=\tau_{1}\tau$,  then we obtain $ h_{\tau'(2)}\cdots h_{\tau'(m)}\in\mathfrak{L}(V), h_{\tau'(2)}\cdots h_{\tau'(m)}\circ h_{\tau'(1)}$ and $ h_{1}\cdots  h_{m}\sim  h_{\tau'(1)}\cdots  h_{\tau'(m)}$ for some $\tau'\in \mathbb  S_{1, \ldots, m}$.

(b)$_{2}$. If all $j_{1}\in\{2, \ldots, b\}, l_{2}\in\{b+1, \ldots, m\}$ such that $\tilde{p}_{ h_{\tau(j_{1})},  h_{j_{2}}}=1$,  then
$ h_{\tau(2)}\cdots  h_{\tau(b)} $ $h_{b+1}\cdots  h_{m}\sim  h_{b+1}\cdots h_{m} h_{\tau(2)}\cdots  h_{\tau(b)}
$. On the other hand,  we know $h_{\tau(1)}\circ h_{b+1}\cdots  h_{m}$ since $ h_{1}\cdots  h_{b}\circ h_{b+1}\cdots  h_{m}$,  then
$ h_{\tau(1)} h_{b+1}\cdots  h_{m}\in\mathfrak L( V)$ by induction hypotheses of (ii). We know $ h_{1}\cdots  h_{m}\sim  h_{\tau(1)}\cdots  h_{\tau(b)} h_{b+1}\cdots  h_{m}\sim  h_{\tau(1)} h_{b+1}\cdots  h_{m} h_{\tau(2)}\cdots h_{\tau(b)}$, and $h_{\tau(1)} h_{b+1}\cdots  h_{m} \circ h_{\tau(2)}\cdots h_{\tau(b)}$ since $h_{\tau(1)}\circ h_{\tau(2)}\cdots h_{\tau(b)}$.
It is proved since
$| h_{\tau(2)}\cdots  h_{\tau(b)}|=b-1<b$ and induction hypotheses.

{\rm (2)}. If  $ h_{1}\cdots  h_{b}= 0$,  then  $ h_{b+1}\cdots  h_{m}\neq0$ and $ h_{b+1}\cdots  h_{m}$ is connected by Lemma \ref {2.3}{\rm (iii)}. Consequently,   $0= h_{1}\cdots  h_{m-1}  \in \mathfrak L(V)$ and $h_{1}\cdots  h_{m-1}\circ h_m.$

(ii) Assume that $ h_{1} h_{2}\cdots  h_{m}\notin\mathfrak L( V)$. Obviously,
$h_{1} h_{2}\cdots  h_{m-1} \not= 0$ and $h_{2} h_{3}\cdots  h_{m} \not=0.$
Set $i _0 :=1$ and   $j_0 := m-1$ when $h_{1}h_{2}\cdots h_{m-1}\in\mathfrak L( V)$;
 $i_0 :=2$ and   $j _0 := m$ when $h_{2}h_{3}\cdots h_{m}\in\mathfrak L( V)$;
$N_0 := \{ i_0, i_0+1, \cdots,  j_0\};  $ $A_0   :=\{1, \ldots, m\}-N_0$,

Now we prove the following $\mathbf{Assertion (k )}$ by induction on $k$,  $1\le k \le m-2.$

$\mathbf{Assertion (k ):}$ There exist $1\le i_k \le  j_k\le m$,  $\tau _k  \in \mathbb S_{N_{k-1}}$ such that the following conditions hold:
$(C_1).$ $0\neq h_{\tau^{k}(i_k)}\cdots h_{\tau^{k}(j_k)}\in\mathfrak L( V);$
$(C_2).$ $\tilde{p}_{h_{r}, h_{t}}=1$ for $\forall \ r\neq t\in A_k$; $(C_3).$
$h_{\tau^{k}(i_k)}\cdots h_{\tau^{k}(j_k)}\circ h_r$ for  $\forall\ r\in A_k$; $(C_4).$
$\tilde{p}_{ h_{\tau^{k}(i_k)}\cdots h_{\tau^{k}(j_k)}   ,  h_{r}}=1$  for  $\forall\ r\in A_k$; $(C_5).$ $h_{\tau^{k}(i_{k-1})}h_{\tau^{k}(i_{k-1}+1)}\cdots h_{\tau^{k}(j_{k-1})} \sim  h_{\tau^{k-1}(i_{k-1})}h_{\tau^{k-1}(i_{k-1}+1)}\cdots h_{\tau^{k-1}(j_{k-1})}  $
 and  $h_{\tau^{k}(1)} \cdots  h_{\tau^{k}(m)} \sim   h_1 h_2 \cdots h_m$; $(C_6).$
 $i_{k-1} \le i_k\le j_{k} \le j_{k-1}$ with  $j_{k-1}- i_{k-1} = j_k -i_k +1$,  where  $N_k := \{\tau^{k} ({i_k}),  \tau^{k}({i_k+1}),  \cdots,  \tau^{k}({j_k})\} $
and  $A_k   :=\{1, \ldots, m\}-N_k$,  $\tau^{k}:=\tau_{k}\tau_{k-1}\cdots \tau_{1}$.

 Step 1. For $k=1, $
 now we construct $i _1$ and $j_1$ as follows.

We know $ 0 \not= h_{i_0}\cdots h_{j_0}\in\mathfrak{L}(V)$, then
 $\tilde{p}_{h_{i_0}\cdots h_{j_0}, h_{r}}=1$  for  $r \in A_1$ by \cite [Lemma 4.12]{WZZ15a} and there exist $\tau_{1}\in \mathbb S_{N_0}$ such that $h_{i_0}\cdots h_{j_0}\sim h_{\tau_{1}(i_0)}\cdots h_{\tau_{1}(j_0)}$ with  (a) $0\neq h_{\tau_{1}(i_0)}\cdots h_{\tau_{1}(j_0-1)}\in\mathfrak{L}(V), $ $i_1 := i_0$ and  $j_1 := j_0-1$
or with (b) $0\neq h_{\tau_{1}(i_0 +1)}\cdots h_{\tau_{1}(j_0)}\in\mathfrak{L}(V), $   $i_1 := i_0 +1$ and $j_1 := j_0$ by induction hypotheses of (iii). Obviously,
$0\neq h_{\tau^{1}(i_1)}\cdots h_{\tau^{1}(j_1)}\in\mathfrak{L}(V)$ and
$h_{\tau_{1}(r)}=h_{r}$ for $r\in A_0$.

Obviously,  $(C_1)$,  $(C_5)$ and $(C_6)$ hold.

$\tilde{p}_{ h_{t},  h_{r}} = 1$ for  any $t \not= r  \in A_1, $ i.e. $(C_2)$ holds.  Indeed,  if $\tilde{p}_{ h_{t},  h_{r}}\neq1$ for $t \in A_1- A_0,  $  $r \in A_0$,  then $\tilde{p}_{ h_{i_1}\cdots  h_{j_1},  h_{r}}\neq1$ since $\tilde{p}_{ h_{i_0}\cdots  h_{j_0},  h_{r}}=1$.
We obtain
$ h_{i_1}\cdots  h_{j_1} h_{t} h_{r}, $  $ h_t h_{i_1}\cdots  h_{j_1} h_{r},  $
$h_r h_{i_1}\cdots  h_{j_1} h_{t} $ $
\in\mathfrak L( V)$ by Lemma \ref {2.1}. However,  $h_{1}\cdots  h_{m}$ is
a quantum equivalent with one among
$ h_{i_1}\cdots  h_{j_1} h_{t} h_{r}, $  $ h_t h_{i_1}\cdots  h_{j_1} h_{r},  $
$h_r h_{i_1}\cdots  h_{j_1} h_{t} $,  which  is a contradiction to
   $h_{1}\cdots  h_{m}\notin\mathfrak L( V)$.

   $(C_3)$ and $(C_4)$ follow from $(C_2)$.

Step 2.  Assumed that $\mathbf{Assertion (k )}$ hold,  we prove that $\mathbf{Assertion (k+1)}$ holds,  $k \leq m-3$.
Let $i:=i_k$ and $j:=j_k$ in this proof for convenience.

We know $0\neq  h_{\tau^{k}(i)}\cdots  h_{\tau^{k}(j)}\in\mathfrak L( V)$,  then there exists
$\tau_{k+1}\in \mathbb  S_{ N_k} $ such that

\noindent $ h_{\tau^{k}(i)}\cdots h_{\tau^{k}(j)}\sim  h_{\tau^{k+1}(i)}\cdots  h_{\tau^{k+1}(j)}$ with ($1^{\circ}$)
$0\neq  h_{\tau^{k+1}(i)}\cdots  h_{\tau^{k+1}(j-1)}\in\mathfrak L( V), $ $i_{k+1} := \tau^{k+1}(i)$  and $j_{k+1}:= \tau^{k+1}(j-1)$,
or with ($2^{\circ}$)
\noindent $0\neq  h_{\tau^{k+1}(i+1)}\cdots  h_{\tau^{k+1}(j)}\in\mathfrak L( V), $ $ i_{k+1} := \tau^{k+1}(i+1)$  and $j_{k+1}:= \tau^{k+1}(j)$ by induction hypotheses of (iii). We obtain
$ h_{1}\cdots  h_{m}\sim  h_{\tau^{k+1}(1)}\cdots  h_{\tau^{k+1}(m)}$.
For convenience,  $ h_{s}':= h_{\tau^{k+1}(s)}$ for all $s\in\{1, \ldots, m\}$. We know $\tau_{k+1}(t)
=t$ for all $t\in A_k $ by the definition of $\tau^{k+1}$.

Obviously,  $(C_1)$,  $(C_5)$ and $(C_6)$ hold.

($1^{\circ}$) If there exists $r\in A_{k+1} $ such that $\tilde{p}_{ h_{\alpha}',  h_{r}}= 1$ for all $\alpha\in\{i, \ldots, j-1\}$,  then
$\tilde{p}_{ h_{j}',   h_{r}}\neq 1$ by $\mathbf{Assertion (k )}$,  we know $\tilde{p}_{ h_{i}'\cdots  h_{j}',   h_{r}}\neq 1$,  it is
a contradiction to $\mathbf{Assertion (k )}$. Thus for all $r\in A_{k+1} $, there exists $\alpha\in\{i, \ldots, j-1\}$ such that $\tilde{p}_{ h_{\alpha}',   h_{r}}\neq 1$ ,  we obtain

\noindent $ h_{1}'\cdots  h_{i-1}' h_{i}'\cdots  h_{j-1}'\in\mathfrak L( V)$ by induction hypotheses.

Set $ \{g_1 , \cdots ,  g_{\beta},  l_1 ,  \cdots , l_{\eta}\} = \{\tau^{k+1}(j+1), \ldots, \tau^{k+1}(m-1)\}$  with $g_1 < \cdots < g_{\beta},  l_1 < \cdots < l_{\eta}$ such that
$\tilde{p}_{ h_{j}',   h_{g_{\lambda}}}\not= 1$ and $\tilde{p}_{ h_{j}',   h_{l_{\xi}}}= 1$ for all $ \lambda  \in\{1, \ldots, \beta\}$,  all \ $\xi\in\{1, \ldots, \eta\}$.
Then $ h_{1}'\cdots  h_{j-1}' h_{l_{1}}'\cdots  h_{l_{r}}'\in\mathfrak L( V)$
for all $r\in\{1, \ldots, \eta\}$ and $ h_{j}' h_{g_{1}}'\cdots h_{g_{t}}'\in\mathfrak L( V)$
for all $t\in\{1, \ldots, \beta\}$ by induction hypotheses.
If $\tilde{p}_{ h_{j}',   h_{m}'}\not= 1$,  then
$\tilde{p}_{ h_{1}'\cdots  h_{j-1}' h_{l_{1}}'\cdots  h_{l_{\eta}}',   h_{m}'}\not= 1$ and
$\tilde{p}_{ h_{j}' h_{g_{1}}'\cdots  h_{g_{\beta}}',   h_{m}'}\not= 1$ by $\mathbf{Assertion (k )}$. Then $ h_{1}'\cdots  h_{j-1}' h_{l_{1}}'\cdots h_{l_{\eta}}' h_{j}' h_{g_{1}}'\cdots  h_{g_{\beta}}' h_{m}'
\in  \mathfrak L( V)$ by Lemma \ref {2.1},
$ h_{1}\cdots  h_{m}\sim  h_{1}'\cdots  h_{m}' \sim   h_{1}'\cdots h_{j-1}' h_{l_{1}}'\cdots  h_{l_{\eta}}' h_{j}' h_{g_{1}}'\cdots  h_{g_{\beta}}' h_{m}'$. It is a contradiction. We obtain $\tilde{p}_{ h_{j}',   h_{m}'}= 1$.
We obtain $ h_{1}'\cdots  h_{j-1}' h_{l_{1}}'\cdots h_{l_{\eta}}' h_{m}'\in\mathfrak{L}(V)$ by induction hypotheses.
If $\beta\geq1$,  then
$\tilde{p}_{ h_{j}' h_{g_{1}}'\cdots  h_{g_{\beta-1}}',  h_{g_{\beta}}'}\neq1$ and
$\tilde{p}_{ h_{1}'\cdots  h_{j-1}' h_{l_{1}}'\cdots h_{l_{\eta}}' h_{m}',  h_{g_{\beta}}'}$
$\neq1$ by $\mathbf{Assertion (k )}$. then
$ h_{1}'\cdots  h_{j-1}' h_{l_{1}}'\cdots  h_{l_{\eta}}' h_{m}' h_{j}' h_{g_{1}}'\cdots  h_{g_{\beta}-1}' h_{g_{\beta}}'\in\mathfrak{L}(V)$ by Lemma \ref {2.1}.
$ h_{1} h_{2}\cdots  h_{m}\sim  h_{1}'\cdots  h_{j-1}' h_{l_{1}}'\cdots h_{l_{\eta}}' h_{m}' h_{j}' h_{g_{1}}'\cdots  h_{g_{\beta}-1}' h_{g_{\beta}}'$. It is contradiction. Then $\beta=0$.
$\eta=m-j-1$,  i.e. $\tilde{p}_{ h_{j}',  h_{r}'}=1$ for all $r\in\{j+1, \ldots, m\}$.
So $\tilde{p}_{ h_{i}'\cdots  h_{j-1}',  h_{r}'}=1$ for all $r\in\{j+1, \ldots, m\}$ by $\mathbf{Assertion (k )}$.
Assumed that there exists $\theta\in\{1, \ldots, i-1\}$ such that $\tilde{p}_{ h_{\theta}',  h_{j}'}\neq1$. If
$\tilde{p}_{ h_{j}',  h_{1}'\cdots  h_{\theta-1}' h_{\theta+1}'\cdots h_{i-1}' h_{i}'\cdots  h_{j-1}'  h_{j+1}'\cdots  h_{m}'}=1$,
then $\tilde{p}_{ h_{j}',  h_{1}'\cdots  h_{j-1}'  h_{j+1}'\cdots  h_{m}'}\neq1$,  $ h_{1}'\cdots  h_{j-1}'  h_{j+1}'\cdots  h_{m}'
\in\mathfrak{L}(V)$ is clear by induction hypotheses. then $ h_{1}'\cdots  h_{j-1}' h_{j+1}'\cdots  h_{m}' h_{j}'$

\noindent $\in\mathfrak{L}(V)$ by \cite [Lemma 4.12]{WZZ15a}. $ h_{1} h_{2}\cdots  h_{m}\sim h_{1}'\cdots  h_{j-1}'  h_{j+1}'\cdots  h_{m}' h_{j}'$. It is contradiction. If
$\tilde{p}_{ h_{j}',  h_{1}'\cdots  h_{\theta-1}' h_{\theta+1}'\cdots h_{i-1}' h_{i}'\cdots  h_{j-1}'  h_{j+1}'\cdots  h_{m}'}\neq1$,
$ h_{1}'\cdots  h_{\theta-1}' h_{\theta+1}'\cdots  h_{i-1}' h_{i}'\cdots h_{j-1}'  h_{j+1}'\cdots  h_{m}'
\in\mathfrak{L}(V)$ is clear by induction hypotheses. Then
$ h_{\theta}' h_{1}'\cdots  h_{\theta-1}' h_{\theta+1}'\cdots  h_{i-1}' h_{i}'\cdots  h_{j-1}'  h_{j+1}'\cdots  h_{m}' h_{j}'$

\noindent $\in\mathfrak{L}(V)$ by Lemma \ref {2.1}.
$ h_{1} h_{2}\cdots  h_{m}\sim  h_{\theta}' h_{1}'\cdots h_{\theta-1}' h_{\theta+1}'\cdots  h_{i-1}' h_{i}'\cdots  h_{j-1}'  h_{j+1}'\cdots  h_{m}' h_{j}'$. It is contradiction. Then
$\tilde{p}_{ h_{\theta}',  h_{j}'}=1$ for all $\theta\in\{1, \ldots, i-1\}$,  which implies that $(C_2)$ holds.

 $(C_3)$ and $(C_4)$ follows from $(C_2)$.

($2^{\circ}$) If there exists $r\in A_{k+1} $ such that $\tilde{p}_{ h_{\alpha}',  h_{r}}= 1$ for all $\alpha\in\{i+1, \ldots, j\}$,  then
$\tilde{p}_{ h_{i}',   h_{r}}\neq 1$ by $\mathbf{Assertion (k )}$,  we know $\tilde{p}_{ h_{i}'\cdots  h_{j}',   h_{r}}\neq 1$,  it is
a contradiction to $\mathbf{Assertion (k )}$. Thus for all $r\in A_{k+1} $, there exists $\alpha\in\{i+1, \ldots, j\}$ such that $\tilde{p}_{ h_{\alpha}',   h_{r}}\neq 1$ ,  we obtain

\noindent $ h_{i+1}'\cdots  h_{j}' h_{j+1}'\cdots  h_{m}'\in\mathfrak L( V)$ by induction hypotheses.

Set $ \{g_1 , \cdots ,  g_{\beta},  l_1 ,  \cdots , l_{\eta}\} = \{\tau^{k+1}(2), \ldots, \tau^{k+1}(i-1)\}$  with $g_1 < \cdots < g_{\beta},  l_1 < \cdots < l_{\eta}$ such that
$\tilde{p}_{ h_{i}',   h_{g_{\lambda }}}\not= 1$ and $\tilde{p}_{ h_{i}',   h_{l_{\xi}}}= 1$ for all $ \lambda \in\{1, \ldots, \beta\}$,  all $\xi\in\{1, \ldots, \eta\}$.
Then $ h_{l_{1}}'\cdots  h_{l_{r}}' h_{i+1}'\cdots  h_{m}'\in\mathfrak L( V)$
for all $r\in\{1, \ldots, \eta\}$ and $ h_{g_{1}}'\cdots  h_{g_{t}}' h_{i}'\in\mathfrak L( V)$
for all $t\in\{1, \ldots, \beta\}$ by induction hypotheses.
If $\tilde{p}_{ h_{i}',   h_{1}'}\not= 1$,  then
$\tilde{p}_{ h_{l_{1}}'\cdots  h_{l_{\eta}}' h_{i+1}'\cdots  h_{m}',   h_{1}'}\not= 1$ and
$\tilde{p}_{ h_{g_{1}}'\cdots  h_{g_{\beta}}' h_{i}',   h_{1}'}\not= 1$ by $\mathbf{Assertion (k )}$. Then $ h_{1}' h_{g_{1}}'\cdots  h_{g_{\beta}}' h_{i}' h_{l_{1}}'\cdots  h_{l_{\eta}}' h_{i+1}'\cdots  h_{m}'
\in  \mathfrak L( V)$ by Lemma \ref {2.1},
$ h_{1}\cdots  h_{m}\sim  h_{1}'\cdots  h_{m}'
\sim  h_{1}' h_{g_{1}}'\cdots  h_{g_{\beta}}' h_{i}'  h_{l_{1}}'\cdots h_{l_{\eta}}' h_{i+1}'\cdots  h_{m}'$. It is a contradiction. We obtain $\tilde{p}_{ h_{i}',   h_{1}'}= 1$.
We obtain $ h_{1}' h_{l_{1}}'\cdots  h_{l_{\eta}}' h_{i+1}'\cdots h_{m}'\in\mathfrak{L}(V)$ by induction hypotheses.
If $\beta\geq1$,  then
$\tilde{p}_{ h_{g_{2}}'\cdots  h_{g_{\beta}}' h_{i}',  h_{g_{1}}'}\neq1$ and
$\tilde{p}_{ h_{1}' h_{l_{1}}'\cdots  h_{l_{\eta}}' h_{i+1}'\cdots h_{m}',  h_{g_{1}}'}$
$\neq1$ by $\mathbf{Assertion (k )}$. then
$ h_{g_{1}}' h_{g_{2}}'\cdots  h_{g_{\beta}}' h_{i}' h_{1}' h_{l_{1}}'\cdots h_{l_{\eta}}' h_{i+1}'\cdots  h_{m}'\in\mathfrak{L}(V)$ by Lemma \ref {2.1}.
$ h_{1} h_{2}\cdots  h_{m}\sim  h_{g_{1}}' h_{g_{2}}'\cdots h_{g_{\beta}}' h_{i}' h_{1}' h_{l_{1}}'\cdots  h_{l_{\eta}}' h_{i+1}'\cdots  h_{m}'$. It is contradiction. Then $\beta=0$.
$\eta=i-2$,  i.e. $\tilde{p}_{ h_{i}',  h_{r}'}=1$ for all $r\in\{1, \ldots, i-1\}$.
So $\tilde{p}_{ h_{i+1}'\cdots  h_{j}',  h_{r}'}=1$ for all $r\in\{1, \ldots, i-1\}$ by $\mathbf{Assertion (k )}$.
Assumed that there exists $\theta\in\{j+1, \ldots, m\}$ such that $\tilde{p}_{ h_{\theta}',  h_{i}'}\neq1$. If
$\tilde{p}_{ h_{i}',  h_{1}'\cdots  h_{i-1}' h_{i+1}'\cdots  h_{j}' h_{j+1}'\cdots h_{\theta-1}' h_{\theta+1}'\cdots  h_{m}'}=1$,
then $\tilde{p}_{ h_{i}',  h_{1}'\cdots  h_{i-1}'  h_{i+1}'\cdots  h_{m}'}\neq1$,  $ h_{1}'\cdots  h_{i-1}'  h_{i+1}'\cdots  h_{m}'
\in\mathfrak{L}(V)$ is clear by induction hypotheses. then $ h_{i}' h_{1}'\cdots h_{i-1}'  h_{i+1}'\cdots  h_{m}'
\in\mathfrak{L}(V)$ by \cite [Lemma 4.12]{WZZ15a}. $ h_{1} h_{2}\cdots  h_{m}\sim h_{i}' h_{1}'\cdots  h_{i-1}'  h_{i+1}'\cdots  h_{m}'$. It is contradiction. If \\
$\tilde{p}_{ h_{i}',  h_{1}'\cdots  h_{i-1}'  h_{i+1}'\cdots  h_{j}' h_{j+1}'\cdots h_{\theta-1}' h_{\theta+1}'\cdots  h_{m}'}\neq1$,
$ h_{1}'\cdots  h_{i-1}'  h_{i+1}'\cdots  h_{j}' h_{j+1}'\cdots h_{\theta-1}' h_{\theta+1}'\cdots  h_{m}'
\in\mathfrak{L}(V)$ is clear by induction hypotheses. Then
$ h_{i}' h_{1}'\cdots  h_{i-1}'  h_{i+1}'\cdots  h_{j}' h_{j+1}'\cdots h_{\theta-1}' h_{\theta+1}'\cdots  h_{m}' h_{\theta}'
$

\noindent $\in\mathfrak{L}(V)$ by Lemma \ref {2.1}.
$ h_{1} h_{2}\cdots  h_{m}\sim  h_{i}' h_{1}'\cdots  h_{i-1}'  h_{i+1}'\cdots h_{j}' h_{j+1}'\cdots  h_{\theta-1}' h_{\theta+1}'\cdots  h_{m}' h_{\theta}'$. It is contradiction. Then
$\tilde{p}_{ h_{\theta}',  h_{i}'}=1$ for all $\theta\in\{j+1, \ldots, m\}$,
which implies that $(C_2)$ holds.

 $(C_3)$ and $(C_4)$ follows from $(C_2)$.

Step 3. In  $\mathbf{Assertion (m-2)}$,  It is a contradiction by $(C_3)$ and $(C_4)$.

(iii) By   Lemma \ref {2.3},  $u$ is connected. which implies that there exist two connected monomials $v$ and $w$
such that $u\sim  vw$. By inductive assumption,   $v$ and $w$ belong to $\mathfrak L(V).$
Consequently,  (iii) holds by (i).

(iv) By Lemma \ref {p2.3},  $u \sim vw$ such that $v$ and $w$ are connected,
as well as,   $v\circ w.$  By inductive assumption,
$v,  w \in  \mathfrak{L}(V). $ It follows from (i) and (ii).  \hfill $\Box$

\begin {Corollary} \label {2.5''} If $ h_{1} h_{2}\cdots  h_{m}\neq0$,  then $ h_{1} h_{2}\cdots  h_{m}\in\mathfrak{L}(V)$
if and only if $\mu ( h_{1} h_{2}\cdots  h_{m})$ is connected.
\end {Corollary}
\noindent {\it Proof.}  It follows from Lemma \ref {2.3} and Theorem \ref {pp2.3} (iv).
\hfill $\Box$

\section {  A basis of Nichols braided Lie algebras }

In this section we obtain a basis of Nichols braided Lie algebra.

\begin {Lemma}  \label{pp2.6} Assume that $V$ is  a  braided vector space   of diagonal type. If $u\neq0$ is homogeneous elements in $\mathfrak B(V)$,  Then $\mu (u)$  is connected (i.e. every monomial of $u$ is connected ) if and only if $u\in \mathfrak L(V)$.
\end {Lemma}

\noindent {\it Proof.}  The necessity follows from  Corollary \ref {2.5''}. The sufficiency. If $u$ is not connected and $u = \sum \limits _{i=1}^r k_i \sigma _i (u_i)$ with $k_i \in F^*$,  where  $u_i$ is a non-zero disconnected monomial and
 $\sigma _i$ is a method of bracket on letters of $u_i$ for $1\le i \le r$. By Lemma \ref {2.3},  $\sigma _i (u_i)=0$ for $1\le i \le r$ and $u=0$,  which is a contradiction.
 $\Box$
\begin {Lemma} \label {p4.1} If $[u] \in D$,  then $u$ is connected and $u \not=0$.
\end {Lemma}
\noindent {\it Proof.}  By \cite [Cor. 1]{Kh99},  $u\not= 0.$ Obviously  $[u] \not=0$ and $[u]\in\mathfrak L( V)$. Considering Lemma \ref {pp2.6}  we complete the proof. $\Box$

\begin {Theorem} \label {2.6'''''} If  $\mathfrak B(V) $ is Nichols algebra of diagonal type with $\dim V\geq2$,  then the set $\{[u_{1}]^{k_1}[u_{2}]^{k_2}\cdots  [u_{s}]^{k_s}\  \mid \  [u_{i}]\in D,  \mid D \mid = s
; 0 \le k_i <  h_{u_i};  1\le i \le s;   u_s<u_{s-1}<\cdots< u_1, \mu ([u_{1}]^{k_1}[u_{2}]^{k_2}\cdots  [u_{s}]^{k_s}) \hbox { is connect }, \sum \limits_{i=1}^{s}k_{i}>0\}$ is a basis of $\mathfrak L (V)$.
\end {Theorem}
\noindent {\it Proof.} It follows from \cite [Th. 1.4.6]{He05},   and Corollary \ref {2.5''} and Lemma \ref {pp2.6}. $\Box$

\section {The dimension of $\mathfrak L(V) $}

In this section we obtain  dimension of Nichols braided Lie algebra $\mathfrak L(V) $
 with finite Cartan type.

Let $V_{i_1,  \cdots,  i_r}$ denote the braided vector subspace generated by $\{x_{i_1}, x_{i_2},  \cdots,   x_{i_r} \}$ of $V$; $D_{i_1,  \cdots,  i_r}$ denote $ \{[u] \mid [u] \hbox { is a hard super-letter of  } \mathfrak B(V_{i_1,  i_2,  \cdots,  i_r})   \}$;

\noindent $ L_{i_1,  \cdots,  i_r} := \{[u_{1}]^{k_1}[u_{2}]^{k_2}\cdots  [u_{s}]^{k_s}\  \mid \  [u_{j}]\in D_{i_1,  \cdots,  i_r}; 0 \le k_j < {\rm ord} (p_{u_j,  u_j});  1\le j \le s; \mid D_{i_1,  \cdots,  i_r} \mid = s;  u_s<u_{s-1}<\cdots< u_1,  \mu ([u_{1}]^{k_1}[u_{2}]^{k_2}\cdots  [u_{s}]^{k_s}) \hbox { is connected }, \sum \limits_{j=1}^{s}k_{j}>0\}$.

\noindent $  B_{i_1,  \cdots,  i_r} := \{[u_{1}]^{k_1}[u_{2}]^{k_2}\cdots  [u_{s}]^{k_s}\  \mid \  [u_{j}]\in D_{i_1,  \cdots,  i_r}; 0 \le k_j < {\rm ord} (p_{u_j,  u_j});  1\le j \le s; \mid D_{i_1,  \cdots,  i_r} \mid = s, \sum \limits_{j=1}^{s}k_{j}>0 \}$; Let $V_{s,  t}$ denote $V_{s,  s+1,   \cdots,  t}$
 in short; similarly we have  $B_{s,  t}$ and  $L_{s, t}$. Let $ B_{i,  j} := \emptyset $ and  $ L_{i,  j} := \emptyset $ when $i >j.$

\begin {Lemma} \label {p4.2} Assumed that $\mathfrak B(V) $ is connected Nichols algebra of diagonal type with $\dim V>2$ and $\Delta(\mathfrak B(V)) $ is an arithmetic root system. The following hold.

{\rm (i)} If pure  generalized Dynkin diagram  is

$\begin{picture}(100,      15)
\put(27,      1){\makebox(0,     0)[t]{$\bullet$}}
\put(60,      1){\makebox(0,      0)[t]{$\bullet$}}
\put(93,     1){\makebox(0,     0)[t]{$\bullet$}}
\put(159,      1){\makebox(0,      0)[t]{$\bullet$}}
\put(192,     1){\makebox(0,      0)[t]{$\bullet$}}
\put(225,     1){\makebox(0,     0)[t]{$\bullet$}}
\put(28,      -1){\line(1,      0){30}}
\put(61,      -1){\line(1,      0){30}}
\put(130,     1){\makebox(0,     0)[t]{$\cdots\cdots\cdots\cdots$}}
\put(160,     -1){\line(1,      0){30}}
\put(193,      -1){\line(1,      0){30}}
\put(22,     -15){1}
\put(58,      -15){2}
\put(91,      -15){3}
\put(157,      -15){n-2}
\put(191,      -15){n-1}
\put(224,      -15){n}
\end{picture}$\\
\\
then
\begin {eqnarray} \label {ep4.2.1} L_{1, n}&=&B_{1, n}-\cup _{i=1}^{n-2}L_{1, i}B_{i+2, n} \\
&=&B_{1, n}-\cup _{i=2}^{n-2}(L_{1, i}- L_{1,  i-1})B_{i+2, n} - L_{1, 1}B_{3, n}\\
&=& B_{1, n}-\cup _{i=1}^{n-3}L_{1, i}(B_{i+2, n}-B_{i+3, n})-L_{1, n-2}B_{n, n}
. \end {eqnarray}
\begin {eqnarray} \label {eppp4.2.2}
 \mid L_{1, n}\mid &=&\mid B_{1, n}\mid -\sum \limits _{i=1}^{n-2} \mid L_{1, i} \mid(\mid B_{i+2, n}\mid-\mid B_{i+3, n}\mid).
\end  {eqnarray}
\begin {eqnarray} \label {epp4.2.2} =\sum \limits _{j=0}^{int(\frac{n-1}{2}) }(-1)^j u_j,
\end {eqnarray}
where
$u _0= \mid B_{1, n}\mid $ and
$ u_j = \sum \limits _{n_1=1}^{n -2}\sum \limits _{n_2=1}^{n_{1} -2 }\cdots \sum \limits _{n_j=1}^{n_{j-1} -2} \mid B_{1, n_j}\mid ( \mid B_{n_j+2, n_{j-1}}\mid
-\mid B_{n_j+3, n_{j-1}}\mid   ) \cdots ( \mid B_{n_2+2, n_1}\mid
-\mid B_{n_2+3, n_1}\mid   )( \mid B_{n_1+2, n}\mid
-\mid B_{n_1+3, n}\mid   )$ for $j>0$.

{\rm (ii)} If pure  generalized Dynkin diagram  is
\\

$\begin{picture}(100,     15)
\put(27,      1){\makebox(0,     0)[t]{$\bullet$}}
\put(75,      1){\makebox(0,      0)[t]{$\bullet$}}
\put(156,     1){\makebox(0,     0)[t]{$\bullet$}}
\put(204,      1){\makebox(0,      0)[t]{$\bullet$}}
\put(245,     -11){\makebox(0,     0)[t]{$\bullet$}}
\put(245,    15){\makebox(0,     0)[t]{$\bullet$}}
\put(28,    -1){\line(1,    0){48}}
\put(115,    1){\makebox(0,     0)[t]{$\cdots\cdots\cdots\cdots$}}
\put(158,     -1){\line(1,     0){48}}
\put(202,    -1){\line(3,     1){42}}
\put(245,      -14){\line(-3,     1){42}}
\put(20,    - 15){1}
\put(70,    - 15){2}
\put(140,     - 15){n - 3}
\put(190,     - 15){n - 2}
\put(250,     8){n - 1}
\put(250,      -18){n}
\end{picture}$\\
\\
then
\begin {eqnarray} \label {ep4.2.3} L_{1, n} &=& B_{1, n} - B_{n-1,  n-1} B_{n,  n} -\cup _{i=1}^{n-3}L_{1, i}B_{i+2, n} \\
&=& B_{1, n}-\cup _{i=2}^{n-3}(L_{1, i} -L_{1,  i-1})B_{i+2, n} -  L_{1,  1}B_{3,  n} - B_{n-1,  n-1} B_{n,  n} \\
&=& B_{1, n}- L_{1,  n-3}B_{n-1,  n} - B_{n-1,  n-1} B_{n,  n}- \cup _{i=1}^{n-3}L_{1, i}(B_{i+2, n} - B_{i+3, n}).
\end {eqnarray}

 \begin {eqnarray} \label {eppp4.2.2'}  \mid L_{1, n}\mid &=&\mid B_{1, n}\mid -\mid B_{n-1, n-1}\mid\mid B_{n, n}\mid-\mid L_{1, n-3}\mid\mid B_{n-1, n}\mid
\nonumber \\
&&
  -\sum \limits _{i=1}^{n-3} \mid L_{1, i} \mid(\mid B_{i+2, n}\mid-\mid B_{i+3, n}\mid). \end  {eqnarray}
where $\mid L_{1, i} \mid $ is obtained by the formula (\ref {epp4.2.2}) when $1\le i \le n-3$.

{\rm (iii)} If pure  generalized Dynkin diagram  is
\\
\\
\\
$\begin{picture}(100,      15)
\put(27,      1){\makebox(0,     0)[t]{$\bullet$}}
\put(60,      1){\makebox(0,      0)[t]{$\bullet$}}
\put(93,     1){\makebox(0,     0)[t]{$\bullet$}}
\put(159,      1){\makebox(0,      0)[t]{$\bullet$}}
\put(192,     1){\makebox(0,      0)[t]{$\bullet$}}
\put(225,     1){\makebox(0,     0)[t]{$\bullet$}}
\put(28,      -1){\line(1,      0){30}}
\put(61,      -1){\line(1,      0){30}}
\put(130,     1){\makebox(0,     0)[t]{$\cdots\cdots\cdots\cdots$}}
\put(160,     -1){\line(1,      0){30}}
\put(193,      -1){\line(1,      0){30}}
\put(22,     -15){1}
\put(58,      -15){2}
\put(91,      -15){3}
\put(157,      -15){n-3}
\put(191,      -15){n-1}
\put(224,      -15){n}

\put(158,    1){\line(0,     1){33}}
\put(164,    28){n-2}
\put(155,    30){$\bullet$}

\end{picture}$\\

then  \begin {eqnarray} \label {ep4.2.5} L_{1, n}&=&B_{1, n}-\cup _{i=1}^{n-4}L_{1, i}B_{i+2, n} - L_{1, n-2}B_{n,  n} - B_{n-2, n-2}
 (B_{n-1, n}-B_{n,  n}) \\
 &=&B_{1, n}-\cup _{i=2}^{n-4}(L_{1, i}-L_{1,  i-1})B_{i+2, n} - L_{1, 1}B_{3,  n} \nonumber  \\
 && -
 (L_{1, n-2} -  L_{1, n-4})B_{n,  n} - B_{n-2, n-2}
 (B_{n-1, n}-B_{n,  n})
 . \end {eqnarray}
\begin {eqnarray} \label {ep4.2.6}
 \mid L_{1, n}\mid &=& \mid B_{1, n}\mid-\sum\limits_{i=2}^{n-4}(\mid  L_{1, i}\mid-\mid L_{1,  i-1}\mid) \mid B_{i+2, n}\mid   -\mid L_{1, 1}\mid \mid B_{3,  n} \mid \nonumber  \\
 && -
(\mid L_{1, n-2}\mid -  \mid L_{1, n-4}\mid)\mid B_{n,  n} \mid- \mid B_{n-2, n-2}\mid(\mid
 B_{n-1, n}\mid- \mid B_{n,  n} \mid). \end {eqnarray}
 where $\mid L_{1, i} \mid $ is obtained by the formula (\ref {epp4.2.2}) when $1\le i \le n-2$.

\end {Lemma}
\noindent {\it Proof.}
 (i) We only determine which element in $B _{1,  n}$ is connected.
It is clear that the left hand of (\ref {ep4.2.1}) $\subseteq $  the right hand of (\ref {ep4.2.1}). If $u\in B _{1,  n}- L _{1,  n}$,   let $ i_u := {\rm min} \{ j \mid x_j \notin \mu (u) \hbox { and there exists } x_i \in \mu (u) \hbox { such that }  1\le i < j \le n\}$. By Lemma \ref {p4.1},  there exist $v \in L_{1,  i_u-1}$ and  $w \in B_{i_u +1,  n }$ such that $u= vw$. Consequently,  the right hand of (\ref {ep4.2.1}) $\subseteq $  the left hand of (\ref {ep4.2.1}). therefore (\ref {ep4.2.1}) holds.

 \begin {eqnarray*}  \mid L_{1, n}\mid &=&\mid B_{1, n}\mid -\sum \limits _{n_1=1}^{n-2} \mid L_{1, n_1} \mid(\mid B_{n_1+2, n}\mid-\mid B_{n_1+3, n}\mid)\\
&=&  \mid B_{1, n}\mid \\
&-& \sum \limits _{n_1=1}^{n-2} (\mid B_{1, n_1}\mid -\sum \limits _{n_2=1}^{n_1-2} \mid L_{1, n_2} \mid(\mid B_{n_2+2, n_1}\mid-\mid B_{n_2+3, n_1}\mid))(\mid B_{n_1+2, n}\mid-\mid B_{n_1+3, n}\mid)\\
&\cdots & \\
&=&\sum \limits _{j=0}^{int(\frac{n-1}{2}) }(-1)^j u_j.
\end {eqnarray*}

Similarly,  we can show (ii) and (iii).  $\Box$

By \cite {Hu78},  $\mid D(A_n) \mid = C_{n+1} ^2, $ $\mid D(B_n) \mid = n^2 = \mid D(C_n) \mid, $ $\mid D(D_n) \mid = n^2-n,\mid D(E_6) \mid =36$,  $\mid D(E_7) \mid =63$,  $\mid D(E_8) \mid =120$,  $\mid D(F_4) \mid =24$,  $\mid D(G_2) \mid =6.$
By Lemma \ref {p4.2} and \cite [Lemma 6.4] {WZZ15a}, we have the following results.

\begin {Example} \label {1}  Let ${\rm ord } (q) :=N$.\\

 (i)
For $A_n$,  $n\ge 1$,
$\begin{picture}(100,      15)
\put(27,      1){\makebox(0,     0)[t]{$\bullet$}}
\put(60,      1){\makebox(0,      0)[t]{$\bullet$}}
\put(93,     1){\makebox(0,     0)[t]{$\bullet$}}
\put(159,      1){\makebox(0,      0)[t]{$\bullet$}}
\put(192,     1){\makebox(0,      0)[t]{$\bullet$}}
\put(28,      -1){\line(1,      0){30}}
\put(61,      -1){\line(1,      0){30}}
\put(130,     1){\makebox(0,     0)[t]{$\cdots\cdots\cdots\cdots$}}
\put(160,     -1){\line(1,      0){30}}
\put(22,     -15){1}
\put(58,      -15){2}
\put(91,      -15){3}
\put(157,      -15){n-1}
\put(191,      -15){n}
\put(22,     10){$q$}
\put(58,      10){$q$}
\put(91,      10){$q$}
\put(157,      10){$q$}
\put(191,      10){$q$}
\put(40,      5){$q^{-1}$}
\put(73,      5){$q^{-1}$}
\put(172,     5){$q^{-1}$}
\put(210,        -1)  {$, q \in F^{*}/\{1\}$. }
\end{picture}$\\
\\ then
$\dim \mathfrak L(V)=\sum \limits _{j=0}^{int(\frac{n-1}{2})} (-1)^j u_j$,  where
$u _0= \mid B_{1,n}\mid $ and
$ u_j = \sum \limits _{n_1=1}^{n -2}\sum \limits _{n_2=1}^{n_{1} -2 }\cdots \sum \limits _{n_j=1}^{n_{j-1} -2} \mid B_{1, n_j}\mid ( \mid B_{n_j+2, n_{j-1}}\mid
-\mid B_{n_j+3, n_{j-1}}\mid   ) \cdots ( \mid B_{n_2+2, n_1}\mid
-\mid B_{n_2+3, n_1}\mid   )( \mid B_{n_1+2, n}\mid
-\mid B_{n_1+3, n}\mid   )$  for $j>0$
 and   $\mid B_{i,  k} \mid  = N ^{C_{k-i +2} ^2}-1$ for
$1\le i \le k \le  n$.
Furthermore,
\begin {eqnarray} \label {epppp4.2.2} && \dim \mathfrak L(V) = \mid L_{1, n} \mid  = N^ {C_{n+1}^2} -1
+
\sum \limits _{j=1}^{int(\frac{n-1}{2})} (-1)^j \sum \limits _{n_1=1}^{n -2}\sum \limits _{n_2=1}^{n_{1} -2 }\cdots \sum \limits _{n_j=1}^{n_{j-1} -2}
( N ^{C_{n_j +1} ^2}-1  )\nonumber \\
&&
(N ^{C_{n_{j-1}- n_j } ^2} - N ^{C_{n_{j-1}- n_j -1} ^2})\cdots (N ^{C_{n_1- n_2 } ^2} - N ^{C_{n_1- n_2 -1} ^2})(N ^{C_{n- n_1 } ^2} - N ^{C_{n- n_1 -1} ^2}).
\end {eqnarray}
\\
(ii) For $B_n$,  $n\ge 2$,
$\begin{picture}(100,      15)
\put(27,      1){\makebox(0,     0)[t]{$\bullet$}}
\put(60,      1){\makebox(0,      0)[t]{$\bullet$}}
\put(93,     1){\makebox(0,     0)[t]{$\bullet$}}
\put(159,      1){\makebox(0,      0)[t]{$\bullet$}}
\put(192,     1){\makebox(0,      0)[t]{$\bullet$}}
\put(225,     1){\makebox(0,     0)[t]{$\bullet$}}
\put(28,      -1){\line(1,      0){33}}
\put(61,      -1){\line(1,      0){30}}
\put(130,     -1){\makebox(0,     0)[t]{$\cdots\cdots\cdots\cdots$}}
\put(160,     -1){\line(1,      0){30}}
\put(193,      -1){\line(1,      0){30}}
\put(22,     -15){1}
\put(58,      -15){2}
\put(91,      -15){3}
\put(157,      -15){n-2}
\put(191,      -15){n-1}
\put(224,      -15){n}
\put(22,     10){$q^{2}$}
\put(58,      10){$q^{2}$}
\put(91,      10){$q^{2}$}
\put(157,      10){$q^{2}$}
\put(191,      10){$q^{2}$}
\put(224,      10){$q$}
\put(40,      5){$q^{-2}$}
\put(73,      5){$q^{-2}$}
\put(172,     5){$q^{-2}$}
\put(205,      5){$q^{-2}$}
\put(235,        -1)  {$, q \in F^{*}/\{1, -1\}$. }
\end{picture}$\\
\\  then
$\dim \mathfrak L(V)=\sum \limits _{j=0}^{int(\frac{n-1}{2})} (-1)^j u_j$,  where
$u _0= \mid B_{1,n}\mid $ and
$ u_j = \sum \limits _{n_1=1}^{n -2}\sum \limits _{n_2=1}^{n_{1} -2 }\cdots \sum \limits _{n_j=1}^{n_{j-1} -2} \mid B_{1, n_j}\mid ( \mid B_{n_j+2, n_{j-1}}\mid
-\mid B_{n_j+3, n_{j-1}}\mid   ) \cdots ( \mid B_{n_2+2, n_1}\mid
-\mid B_{n_2+3, n_1}\mid   )( \mid B_{n_1+2, n}\mid
-\mid B_{n_1+3, n}\mid   )$ for $j>0$;   $\mid B_{i,  n} \mid  = N ^{(n-i +1)^2}-1$ for
$1\le i < n$ and $\mid B_{i,  k} \mid  = N ^{C_{k-i +2} ^2}-1$ for
$1\le i \le k <  n$,  when  N is  odd;
  $\mid B_{i,  n} \mid  = (\frac {N} {2}) ^{(n-i +1)^2-n +i-1}   N ^{n -i+1}-1$ for
$1\le i < n$ and $\mid B_{i,  k} \mid  = (\frac {N} {2}) ^{C_{k-i +2} ^2}-1$ for
$1\le i\le k < n$,  when  N is  even. Furthermore,
\begin {eqnarray} \label {epppp4.2.2'} && \dim \mathfrak L(V) = \mid L_{1, n} \mid  =N ^{n^2}-1
+
\sum \limits _{j=1}^{int(\frac{n-1}{2})} (-1)^j \sum \limits _{n_1=1}^{n -2}\sum \limits _{n_2=1}^{n_{1} -2 }\cdots \sum \limits _{n_j=1}^{n_{j-1} -2}
( N ^{C_{n_j +1} ^2}-1  )\nonumber \\
&&
(N ^{C_{n_{j-1}- n_j } ^2} - N ^{C_{n_{j-1}- n_j -1} ^2})\cdots (N ^{C_{n_1- n_2 } ^2} - N ^{C_{n_1- n_2 -1} ^2}) \nonumber\\
 && (N ^{(n- n_1-1)^ 2 - n +n_1 +1 } - N ^{(n- n_1-2)^ 2 - n +n_1 +2 } ),
\end {eqnarray}
 when  N is  odd;
\noindent \begin {eqnarray} \label {epppp4.2.2'} && \dim \mathfrak L(V) = \mid L_{1, n} \mid  = (\frac {N} {2})^{n^2-n }  N ^n-1
+
\sum \limits _{j=1}^{int(\frac{n-1}{2})} (-1)^j \sum \limits _{n_1=1}^{n -2}\sum \limits _{n_2=1}^{n_{1} -2 }\cdots \sum \limits _{n_j=1}^{n_{j-1} -2}
( (\frac {N} {2}) ^{C_{n_j +1} ^2}-1  )\nonumber \\
&&
((\frac {N} {2}) ^{C_{n_{j-1}- n_j } ^2} - (\frac {N} {2}) ^{C_{n_{j-1}- n_j -1} ^2})\cdots ((\frac {N} {2}) ^{C_{n_1- n_2 } ^2} - (\frac {N} {2}) ^{C_{n_1- n_2 -1} ^2})\nonumber\\
 &&   ((\frac {N} {2}) ^{(n-n_{1} -1)^2-n +n_{1}+1}   N ^{n -n_{1}-1} - (\frac {N} {2}) ^{(n-n_{1} -2)^2-n +n_{1}+2}   N ^{n -n_{1}-2}),
\end {eqnarray}
 when  N is even.
\\
(iii) For $C_n$,  $n>2$,
$\begin{picture}(100,      15)
\put(27,      1){\makebox(0,     0)[t]{$\bullet$}}
\put(60,      1){\makebox(0,      0)[t]{$\bullet$}}
\put(93,     1){\makebox(0,     0)[t]{$\bullet$}}
\put(159,      1){\makebox(0,      0)[t]{$\bullet$}}
\put(192,     1){\makebox(0,      0)[t]{$\bullet$}}
\put(225,     1){\makebox(0,     0)[t]{$\bullet$}}
\put(28,      -1){\line(1,      0){33}}
\put(61,      -1){\line(1,      0){30}}
\put(130,     -1){\makebox(0,     0)[t]{$\cdots\cdots\cdots\cdots$}}
\put(160,     -1){\line(1,      0){30}}
\put(193,      -1){\line(1,      0){30}}
\put(22,     -15){1}
\put(58,      -15){2}
\put(91,      -15){3}
\put(157,      -15){n-2}
\put(191,      -15){n-1}
\put(224,      -15){n}
\put(22,     10){$q$}
\put(58,      10){$q$}
\put(91,      10){$q$}
\put(157,      10){$q$}
\put(191,      10){$q$}
\put(224,      10){$q^2$}
\put(40,      5){$q^{-1}$}
\put(73,      5){$q^{-1}$}
\put(172,     5){$q^{-1}$}
\put(205,      5){$q^{-2}$}
\put(235,        -1)  {$, q \in F^{*}/\{1, -1\}$. }
\end{picture}$\\
\\  then
$\dim \mathfrak L(V)=\sum \limits _{j=0}^{int(\frac{n-1}{2})} (-1)^j u_j$,  where
$u _0= \mid B_{1,n}\mid  $ and
$ u_j = \sum \limits _{n_1=1}^{n -2}\sum \limits _{n_2=1}^{n_{1} -2 }\cdots \sum \limits _{n_j=1}^{n_{j-1} -2} \mid B_{1, n_j}\mid ( \mid B_{n_j+2, n_{j-1}}\mid
-\mid B_{n_j+3, n_{j-1}}\mid   ) \cdots ( \mid B_{n_2+2, n_1}\mid
-\mid B_{n_2+3, n_1}\mid   )( \mid B_{n_1+2, n}\mid
-\mid B_{n_1+3, n}\mid   )$ for $j>0$;   $\mid B_{i,  n} \mid  = N ^{(n-i +1)^2}-1$ for
$1\le i < n$ and $\mid B_{i,  k} \mid  = N ^{C_{k-i +2} ^2}-1$ for
$1\le i\le k < n$,  when  N is  odd;
  $\mid B_{i,  n} \mid  = N ^{(n-i +1)^2-n +i-1}   (\frac {N} {2}) ^{n -i +1}-1$ for
$1\le i < n$ and $\mid B_{i,  k} \mid  = N ^{C_{k-i +2} ^2}-1$ for
$1\le i\le k < n$,  when  N is  even.
Furthermore,
\begin {eqnarray} \label {epppp4.2.2''} && \dim \mathfrak L(V) = \mid L_{1, n} \mid  =N ^{n^2 }   -1
+
\sum \limits _{j=1}^{int(\frac{n-1}{2})} (-1)^j \sum \limits _{n_1=1}^{n -2}\sum \limits _{n_2=1}^{n_{1} -2 }\cdots \sum \limits _{n_j=1}^{n_{j-1} -2}
( N ^{C_{n_j +1} ^2}-1  )\nonumber \\
&&
(N ^{C_{n_{j-1}- n_j } ^2} - N ^{C_{n_{j-1}- n_j -1} ^2})\cdots (N ^{C_{n_1- n_2 } ^2} - N ^{C_{n_1- n_2 -1} ^2})
\nonumber \\
&&(N ^{(n- n_1-1)^ 2 - n +n_1 +1 } - N ^{(n- n_1-2)^ 2 - n +n_1 +2 } ),
\end {eqnarray}
 when  N is  odd;
\begin {eqnarray} \label {epppp4.2.2''''} && \dim \mathfrak L(V) = \mid L_{1, n} \mid  =(\frac {N} {2}) ^{n }   N ^{n^2-n }-1
+
\sum \limits _{j=1}^{int(\frac{n-1}{2})} (-1)^j \sum \limits _{n_1=1}^{n -2}\sum \limits _{n_2=1}^{n_{1} -2 }\cdots \sum \limits _{n_j=1}^{n_{j-1} -2}
( {N}  ^{C_{n_j +1} ^2}-1  )\nonumber \\
&&
(N ^{C_{n_{j-1}- n_j } ^2} - N ^{C_{n_{j-1}- n_j -1} ^2})\cdots (N ^{C_{n_1- n_2 } ^2} - N ^{C_{n_1- n_2 -1} ^2})\nonumber \\
&&  (N ^{(n-n_{1} -1)^2-n +n_{1}+1}   (\frac {N} {2}) ^{n -n_{1} -1} - N ^{(n-n_{1} -2)^2-n +n_{1}+2}   (\frac {N} {2}) ^{n -n_{1} -2}),
\end {eqnarray}
 when  N is even.
\\

(iv) For  $D_n$,  $n>3$,
$\begin{picture}(100,     15)
\put(12,      1){\makebox(0,     0)[t]{$\bullet$}}
\put(45,      1){\makebox(0,      0)[t]{$\bullet$}}
\put(111,     1){\makebox(0,     0)[t]{$\bullet$}}
\put(144,      1){\makebox(0,      0)[t]{$\bullet$}}
\put(170,     -11){\makebox(0,     0)[t]{$\bullet$}}
\put(170,    15){\makebox(0,     0)[t]{$\bullet$}}
\put(13,    -1){\line(1,    0){33}}
\put(80,    1){\makebox(0,     0)[t]{$\cdots\cdots\cdots\cdots$}}
\put(113,     -1){\line(1,     0){33}}
\put(142,    -1){\line(2,     1){27}}
\put(170,      -14){\line(-2,     1){27}}
\put(5,    - 15){1}
\put(40,    - 15){2}
\put(95,     - 15){n - 3}
\put(130,     - 15){n - 2}
\put(175,     8){n - 1}
\put(175,      -18){n}
\put(10,      10){$q$}
\put(45,      10){$q$}
\put(110,      10){$q$}
\put(143,      10){$q$}
\put(170,      -8){$q$}
\put(170,      18){$q$}
\put(25,      5){$q^{-1}$}
\put(123,      5){$q^{-1}$}
\put(150,     -4){$q^{-1}$}
\put(150,      10){$q^{-1}$}
\put(195,        -1)  {$, q \in F^{*}/\{1\}$. }
\end{picture}$\\
\\
then $ \dim \mathfrak L(V)=\mid B_{1, n}\mid -\mid B_{n-1, n-1}\mid\mid B_{n, n}\mid-\sum \limits _{i=1}^{n-3} \mid L_{1, i} \mid(\mid B_{i+2, n}\mid-\mid B_{i+3, n}\mid)-\mid L_{1, n-3}\mid\mid B_{n-1, n}\mid
$,
where
$\mid B_{ k ,n}\mid = N^{({n-k+1})^2-n +k -1}-1$ for $n-k +1 >2$,
$\mid B_{ n-1, n}\mid = N^3-1$, $ \mid B_{ n, n}\mid = \mid B_{ n-1, n-1}\mid = N-1$; $\mid L_{1, t} \mid $ is obtained by the formula (\ref {epppp4.2.2}) when $1\le t \le n-3$.
Furthermore,
\begin {eqnarray} \label {epppp4.2.2'''''} && \dim \mathfrak L(V) = \mid L_{1, n} \mid  =
N^{n^2-n} - (N-1)^2
-\sum \limits _{i=1}^{n-3} \{N^ {C_{i+1}^2} -1
+
\sum \limits _{j=1}^{int(\frac{i-1}{2})} (-1)^j \sum \limits _{n_1=1}^{n -2}\sum \limits _{n_2=1}^{n_{1} -2 }\cdots \nonumber \\
&&
\cdots\sum \limits _{n_j=1}^{n_{j-1} -2}( N ^{C_{n_j +1} ^2}-1  )
(N ^{C_{n_{j-1}- n_j } ^2} - N ^{C_{n_{j-1}- n_j -1} ^2})\cdots (N ^{C_{n_1- n_2 } ^2} - N ^{C_{n_1- n_2 -1} ^2})\nonumber \\
&&(N ^{C_{i- n_1 } ^2} - N ^{C_{i- n_1 -1} ^2})\}(N^{({n-i-1})^2-n +i +1}-N^{({n-i-2})^2-n +i +2})-\{N^ {C_{n-2}^2} -1
\nonumber \\
&&
+\sum \limits _{j=1}^{int(\frac{n-4}{2})} (-1)^j \sum \limits _{n_1=1}^{n -5}\sum \limits _{n_2=1}^{n_{1} -2 }\cdots \sum \limits _{n_j=1}^{n_{j-1} -2}( N ^{C_{n_j +1} ^2}-1  )
(N ^{C_{n_{j-1}- n_j } ^2} - N ^{C_{n_{j-1}- n_j -1} ^2})\cdots \nonumber \\
&&
\cdots(N ^{C_{n_1- n_2 } ^2} - N ^{C_{n_1- n_2 -1} ^2})(N ^{C_{n-3- n_1 } ^2} - N ^{C_{n-4- n_1 } ^2})\}(N^2-1).
\end {eqnarray}
\\

 (v) For $E_6$,
$\begin{picture}(100,    15)
\put(47,    1){\makebox(0,   0)[t]{$\bullet$}}
\put(87,    1){\makebox(0,0)[t]{$\bullet$}}
\put(127,   38){\makebox(0,0)[t]{$\bullet$}}
\put(127,    1){\makebox(0,    0)[t]{$\bullet$}}
\put(167,   1){\makebox(0,    0)[t]{$\bullet$}}
\put(207,   1){\makebox(0,   0)[t]{$\bullet$}}
\put(48,   -1){\line(1,    0){37}}
\put(88,  -1){\line(1,    0){37}}
\put(127,   1){\line(0,    1){37}}
\put(128,  -1){\line(1,    0){37}}
\put(168,  -1){\line(1,    0){37}}
\put(45,   -15){1}
\put(85,   -15){2}
\put(125,   -15){3}
\put(119,   26){4}
\put(166,   -15){5}
\put(206,   -15){6}
\put(49,    8){\makebox(0,   0)[t]{$q$}}
\put(89,    8){\makebox(0,0)[t]{$q$}}
\put(133,   40){\makebox(0,0)[t]{$q$}}
\put(129,    8){\makebox(0,    0)[t]{$q$}}
\put(169,   8){\makebox(0,    0)[t]{$q$}}
\put(209,   8){\makebox(0,   0)[t]{$q$}}
\put(60,   5){$q^{-1}$}
\put(101,   5){$q^{-1}$}
\put(128,   15){$q^{-1}$}
\put(143,   5){$q^{-1}$}
\put(180,   5){$q^{-1}$}
\put(227,       -1)  {$,q \in F^{*}/\{1\}$. }
\end{picture}$\\
\\then
\noindent \begin {eqnarray*}
  \dim \mathfrak L(V) &=& \mid B_{1, 6}\mid-\sum\limits_{i=2}^{6-4}(\mid  L_{1, i}\mid-\mid L_{1,  i-1}\mid) \mid B_{i+2, 6}\mid   -\mid L_{1, 1}\mid \mid B_{3,  6} \mid \nonumber  \\
 && -
(\mid L_{1, 6-2}\mid -  \mid L_{1, 6-4}\mid)\mid B_{6,  6} \mid- \mid B_{6-2, 6-2}\mid(\mid
 B_{6-1, 6}\mid- \mid B_{6,  6} \mid) \\
 &=&N ^{36}-1 - (\mid  L_{1,2}\mid-\mid L_{1, 1}\mid )  ( NN^3-1)   - (N-1) (N ^{C_5^2}-1) \\
 &-&
(\mid L_{1,4} \mid- \mid L_{1,2} \mid) (N-1)- (N-1) ( N ^3 - N ),
 \end {eqnarray*}
 where $\mid L_{1, t} \mid $ is obtained by the formula (\ref {epppp4.2.2}) when $1\le t \le 4$.  \\
  \\
  \\

(vi) For $E_7$,
$\begin{picture}(100,    15)
\put(37,    1){\makebox(0,    0)[t]{$\bullet$}}
\put(77,    1){\makebox(0,   0)[t]{$\bullet$}}
\put(117,    1){\makebox(0,0)[t]{$\bullet$}}
\put(157,   38){\makebox(0,0)[t]{$\bullet$}}
\put(157,    1){\makebox(0,    0)[t]{$\bullet$}}
\put(197,   1){\makebox(0,    0)[t]{$\bullet$}}
\put(237,   1){\makebox(0,   0)[t]{$\bullet$}}
\put(38,   -1){\line(1,    0){37}}
\put(78,   -1){\line(1,    0){37}}
\put(118,  -1){\line(1,    0){37}}
\put(157,   1){\line(0,    1){37}}
\put(158,  -1){\line(1,    0){37}}
\put(198,  -1){\line(1,    0){37}}
\put(35,   -15){1}
\put(75,   -15){2}
\put(115,   -15){3}
\put(155,   -15){4}
\put(149,   26){5}
\put(196,   -15){6}
\put(236,   -15){7}
\put(39,    8){\makebox(0,    0)[t]{$q$}}
\put(79,    8){\makebox(0,   0)[t]{$q$}}
\put(119,    8){\makebox(0,0)[t]{$q$}}
\put(163,   40){\makebox(0,0)[t]{$q$}}
\put(159,    8){\makebox(0,    0)[t]{$q$}}
\put(199,   8){\makebox(0,    0)[t]{$q$}}
\put(239,   8){\makebox(0,   0)[t]{$q$}}
\put(50,   5){$q^{-1}$}
\put(90,   5){$q^{-1}$}
\put(131,   5){$q^{-1}$}
\put(158,   15){$q^{-1}$}
\put(173,   5){$q^{-1}$}
\put(210,   5){$q^{-1}$}
\put(257,       -1)  {$,q \in F^{*}/\{1\}$. }
\end{picture}$\\
\\then
\noindent \begin {eqnarray*}
  \dim \mathfrak L(V) &=& \mid B_{1,7}\mid-\sum\limits_{i=2}^{7-4}(\mid  L_{1,i}\mid- \mid L_{1, i-1}\mid) \mid B_{i+2,7}\mid   -\mid L_{1,1}\mid \mid B_{3, 7} \mid \nonumber  \\
 && -
(\mid L_{1,7-2}\mid - \mid L_{1,7-4}\mid)\mid B_{7, 7} \mid- \mid B_{7-2,7-2}\mid
 (\mid B_{7-1,7}\mid - \mid B_{7, 7} \mid)\\
  &=& \mid B_{1,7}\mid-(\mid  L_{1,2}\mid-\mid L_{1, 1}\mid) \mid B_{3,7}\mid -(\mid  L_{1,3}\mid-\mid L_{1, 2}\mid) \mid B_{5,7}\mid
    -\mid L_{1,1}\mid \mid B_{3, 7} \mid \nonumber  \\
 && -
(\mid L_{1,7-2}\mid - \mid L_{1,7-4}\mid)\mid B_{7, 7} \mid- \mid B_{7-2,7-2}\mid
 (\mid B_{7-1,7}\mid- \mid B_{7, 7} \mid)\\
  &=& N^{63}-1
 -  ( \mid  L_{1,2}\mid-\mid L_{1, 1}\mid)  (N^{ 10}-1) -  (\mid  L_{1,3}\mid- \mid L_{1, 2}\mid)   (NN^{ 3}-1)
 \\
 &&   -
    ( N-1)( N^ {20}-1) -(\mid L_{1,5}\mid -\mid  L_{1,3} \mid) (N-1)-  (N-1)( N ^{3}
- N), \\
 \end {eqnarray*}
 where $\mid L_{1, t} \mid $ is obtained by the formula (\ref {epppp4.2.2}) when $1\le t \le 5$.\\
 \\

(vii) For $E_8$,
$\begin{picture}(100,    15)
\put(27,    1){\makebox(0,    0)[t]{$\bullet$}}
\put(67,    1){\makebox(0,    0)[t]{$\bullet$}}
\put(107,    1){\makebox(0,   0)[t]{$\bullet$}}
\put(147,    1){\makebox(0,0)[t]{$\bullet$}}
\put(187,   38){\makebox(0,0)[t]{$\bullet$}}
\put(187,    1){\makebox(0,    0)[t]{$\bullet$}}
\put(227,   1){\makebox(0,    0)[t]{$\bullet$}}
\put(267,   1){\makebox(0,   0)[t]{$\bullet$}}
\put(28,   -1){\line(1,    0){37}}
\put(68,   -1){\line(1,    0){37}}
\put(108,   -1){\line(1,    0){37}}
\put(148,  -1){\line(1,    0){37}}
\put(187,   1){\line(0,    1){37}}
\put(188,  -1){\line(1,    0){37}}
\put(228,  -1){\line(1,    0){37}}
\put(22,  -15){1}
\put(65,   -15){2}
\put(105,   -15){3}
\put(145,   -15){4}
\put(185,   -15){5}
\put(179,   26){6}
\put(226,   -15){7}
\put(266,   -15){8}
\put(29,    8){\makebox(0,    0)[t]{$q$}}
\put(69,    8){\makebox(0,    0)[t]{$q$}}
\put(109,    8){\makebox(0,   0)[t]{$q$}}
\put(149,    8){\makebox(0,0)[t]{$q$}}
\put(193,   40){\makebox(0,0)[t]{$q$}}
\put(189,    8){\makebox(0,    0)[t]{$q$}}
\put(229,   8){\makebox(0,    0)[t]{$q$}}
\put(269,   8){\makebox(0,   0)[t]{$q$}}
\put(40,  5){$q^{-1}$}
\put(80,   5){$q^{-1}$}
\put(120,   5){$q^{-1}$}
\put(161,   5){$q^{-1}$}
\put(188,   15){$q^{-1}$}
\put(203,   5){$q^{-1}$}
\put(240,   5){$q^{-1}$}
\put(287,       -1)  {$,q \in F^{*}/\{1\}$. }
\end{picture}$\\
\\then
\begin {eqnarray*}
  \dim \mathfrak L(V) &=& \mid B_{1,8}\mid-\sum\limits_{i=2}^{8-4}(\mid  L_{1,i}\mid- \mid L_{1, i-1}\mid) \mid B_{i+2,8}\mid   -\mid L_{1,1}\mid\mid B_{3, 8} \mid \nonumber  \\
 && -
(\mid L_{1,8-2}\mid - \mid L_{1,8-4}\mid)\mid B_{8, 8} \mid- \mid B_{8-2,8-2}\mid(\mid
 B_{8-1,8}\mid-\mid B_{8, 8} \mid)\\
 &=&   \mid B_{1,8}\mid  -
(\mid  L_{1,2}\mid- \mid L_{1, 1}\mid )\mid B_{4,8}\mid
  -(\mid  L_{1,3}\mid-\mid L_{1, 2}\mid) \mid B_{5,8}\mid\\
   &&-
   (\mid  L_{1,4}\mid- \mid L_{1, 3}\mid )\mid B_{6,8}\mid
   -\mid L_{1,1}\mid \mid B_{3, 8} \mid
  -(\mid L_{1,6}\mid - \mid L_{1,4}\mid)\mid B_{8, 8} \mid\\
  &&-
   \mid B_{6,6}\mid(\mid B_{7,8}\mid- \mid B_{8, 8}\mid)
\\
 &=&  N  ^{120}-1  -
(\mid  L_{1,2}\mid- \mid L_{1, 1}\mid) ( N^ {20}-1)
  -(\mid  L_{1,3}\mid-\mid L_{1, 2}\mid )  ( N^ {C_5^2}-1)    \\
  &&-
  (\mid  L_{1,4}\mid-\mid L_{1, 3}\mid)  (N N ^3-1)
   - (N-1)(N^{36}-1)  \\
 && -
(\mid L_{1,6}\mid -\mid  L_{1,4} \mid )(N-1)-   (N-1)( N^3 - N )
 , \end {eqnarray*}
 where $\mid L_{1, t} \mid $ is obtained by the formula (\ref {epppp4.2.2}) when $1\le t \le 6$.\\
 \\

(viii) For $F_4$.
$\begin{picture}(100,      15)
\put(27,      1){\makebox(0,     0)[t]{$\bullet$}}
\put(60,      1){\makebox(0,      0)[t]{$\bullet$}}
\put(93,     1){\makebox(0,     0)[t]{$\bullet$}}
\put(126,      1){\makebox(0,    0)[t]{$\bullet$}}
\put(28,      -1){\line(1,      0){33}}
\put(61,      -1){\line(1,      0){30}}
\put(94,     -1){\line(1,      0){30}}
\put(22,     -15){1}
\put(58,      -15){2}
\put(91,      -15){3}
\put(124,     -15){4}
\put(22,     10){$q^2$}
\put(58,      10){$q^2$}
\put(91,      10){$q$}
\put(124,      10){$q$}
\put(40,      5){$q^{-2}$}
\put(73,      5){$q^{-2}$}
\put(106,     5){$q^{-1}$}
\put(145,        -1)  {$, q \in F^{*}/\{1, -1\}$. }
\end{picture}$\\
\\
 then
$\dim \mathfrak L(V)=\sum \limits _{j=0}^{1} (-1)^j u_j$,  where $u_0 =  \mid B_{1, 4}\mid $ and
$ u_1 = \sum \limits _{n_1=1}^{2} ( \mid B_{1, n_1}\mid ) ( \mid B_{n_1+2, 4}\mid
-\mid B_{n_1+3, 4}\mid   )$;  $\mid B_{1, 4}\mid  = N ^{24}-1$,  $\mid B_{1, 1}\mid = N -1$,
 $\mid B_{1, 2}\mid  = N ^{C_3^2}-1 = N^3-1$,  $\mid B_{3, 4 }\mid  = N ^{C_3^2}-1 = N^3-1$ and   $\mid B_{4, 4 }\mid  = N -1 $ when  $N$ is  odd;     $\mid B_{1, 4}\mid  = N ^{12} (\frac {N} {2}) ^{12}-1$,  $\mid B_{1, 1}\mid = \frac {N} {2} -1$,
 $\mid B_{1, 2}\mid  = (\frac {N} {2}) ^{C_3^2}-1 = (\frac {N} {2})^3-1$,  $\mid B_{3, 4 }\mid  = N ^{C_3^2}-1 = N^3-1$ and   $\mid B_{4, 4 }\mid  = N -1 $ when  $N$ is  even.
Furthermore
\begin {eqnarray} \label {epppp4.2.2''''''}  \dim \mathfrak L(V)   = N^ {24} -1 -( N -1)(N ^{3}-N ) - (N^3-1) (N-1),
\end {eqnarray} when $N$ is odd;
\begin {eqnarray} \label {epppp4.2.2.1}  \dim \mathfrak L(V)   = (\frac {N}{2})^ {12} N^{12} -1 -( \frac {N}{2} -1)(N ^{3}-N ) - ((\frac {N}{2})^ {3}-1) (N-1),
\end {eqnarray} when $N$ is even.

(ix) For $G_2$,    $\begin{picture}(100,      15)
\put(27,      1){\makebox(0,     0)[t]{$\bullet$}}
\put(60,      1){\makebox(0,      0)[t]{$\bullet$}}
\put(28,      -1){\line(1,      0){33}}
\put(22,     -15){1}
\put(58,      -15){2}
\put(22,     10){$q$}
\put(58,      10){$q^3$}
\put(40,      5){$q^{-3}$}
  \ \ \ \ \ \ \ \ \ \ \ \ \ \ \ \ \ \ \ {$,q \in F^{*}/\{1, -1\}, q^3 \not= 1$. }
\end{picture}$\\
\\

 then
$\dim \mathfrak L(V)=N^ 6-1$ when $3 \nmid N$; $\dim \mathfrak L(V)=(\frac {N}{3})^3 N^3-1$ when $3 \mid N$.
 \end {Example}

\section { Non-zero monomials in  Nichols algebras }

In this section we find some non-zero monomials in  Nichols algebras.

 Let $(s)_q := 1 + q + q^2 + \cdots + q ^{s-1}$ and $ (s) _q ! := (1)_q (2)_q \cdots (s)_q.$
\begin{Lemma}\label{2.0} Assume $h_i \in \{x_1,  \cdots,  x_n\}$ for $1\le i \le m$ and $u _j$ is 1 or a monomial with $x_k \notin \mu (u_j)$ for $ 1\le j \le l+1$. Let $q:= p_{x_k, x_k} ^{-1}$ for convenience.

{\rm (i)}   \begin {eqnarray} \label {ep201} &&<y_{k}, u_{1}x_{k}u_{2}x_{k}\cdots u_{l}x_{k}u_{l+1}>
 \nonumber \\
&=& \sum \limits _{j=1}^{l}  q^{j-1} p _{k, u_1u_2\cdots u_j} ^{-1} u_1 x_k u_2 x_k \cdots x_k
(u_j u_{j+1}) x_k \cdots x_k u_{l+1}. \end {eqnarray}

{\rm (ii)} \begin {eqnarray} \label {ep202} <y_{k} ^l, u_{1}x_{k}u_{2}x_{k}\cdots u_{l}x_{k}u_{l+1}>
&=&  (l)_q!p_{k, u_{1}}^{-1}p_{k, u_{1}u_{2}}^{-1}\cdots p_{k, u_{1}u_{2}\cdots u_{l}}^{-1}
u_{1}u_{2}\cdots u_{l+1}
. \end {eqnarray}

{\rm (iii)} If ${\rm ord } (p_{h_i,h_i} ) >\mid {\rm deg}_ { h_i}( h_{1}\cdots  h_{m})\mid$ for all $i \in\{1, \ldots, m\}$,
then $ h_{1}\cdots  h_{m}\neq0$.
\end {Lemma}
\noindent {\it Proof.} {\rm (i)}  It can be obtained by induction on $l.$

 {\rm (ii)}
 We show this by induction on $l$. If $l=1$,  $<y_{k}, u_{1}x_{k}u_{2}>=p_{k, u_{1}}^{-1}
u_{1}u_{2}$. Assume $l>1.$  See that
\begin {eqnarray*}
 &&<y_{k}^{l}, u_{1}x_{k}u_{2}x_{k}\cdots u_{l}x_{k}u_{l+1}> \\
    &=&
<y_{k}^{l-1}, <y_{k}, u_{1}x_{k}u_{2}x_{k}\cdots u_{l}x_{k}u_{l+1}>> \\
&=&<y_k ^{l-1},  \sum \limits _{j=1}^{l}  q^{j-1} p _{k, u_1u_2\cdots u_j} ^{-1} u_1 x_k u_2 x_k \cdots x_k
(u_j u_{j+1}) x_k \cdots x_k u_{l+1}  > \\
&=&  \sum \limits _{j=1}^{l}  q^{j-1} p _{k, u_1u_2\cdots u_j} ^{-1}  < y_k ^{l-1},  u_1 x_k u_2 x_k \cdots x_k
(u_j u_{j+1}) x_k \cdots x_k u_{l+1}> \\
&=&  (l)_q!p_{k, u_{1}}^{-1}p_{k, u_{1}u_{2}}^{-1}\cdots p_{k, u_{1}u_{2}\cdots u_{l}}^{-1}
u_{1}u_{2}\cdots u_{l+1}.
\end {eqnarray*}

{\rm (iii)} We show this by induction on $t:= \mid \mu (h_1h_2\cdots h_m)\mid $. If $t=1$,  we obtain $ h_{1}\cdots  h_{m}\neq0$ by \cite [Lemma 1.3.3 {\rm (i)}]{He05}.
Assume $t >1$ and $ h_1h_2\cdots h_m = u_{1}x_{k}u_{2}x_{k}\cdots u_{l}x_{k}u_{l+1}$
 with $\mid \mu (u_{1}u_{2}\cdots u_{l+1})\mid $ $ = t-1$. Thus  $u_{1}u_{2}\cdots u_{l+1}\neq0$ by induction hypotheses,
$(1+p_{kk}^{-1})(1+p_{kk}^{-1}+p_{kk}^{-2})\cdots
(1+p_{kk}^{-1}+\cdots+p_{kk}^{-l+1})\neq0$ since ${\rm ord } (p_{x_{k}, x_{k}})>\mid {\rm deg}_ {x_k}( h_{1}\cdots  h_{m})\mid=l$.
Hence $u_{1}x_{k}u_{2}x_{k}\cdots u_{l}x_{k}u_{l+1}\neq0$,  completing the proof.
 \hfill $\Box$

\begin {Corollary}  \label{2.3'} If $u$ and $ v$ are monomials  in $\mathfrak B(V)$,  $\mu (u)\cap \mu (v)=\emptyset$,  $p_{ij} p_{ji}=1$ for any $x_{i} \in \mu (u)$,  $x_{j} \in \mu (v)$,  then the following conditions are equivalent:
{\rm (i)} $u\neq0, v\neq0$.
{\rm (ii)} $uv\neq0$.
{\rm (iii)} $uv\notin \mathfrak L(V)$.
\end {Corollary}
\noindent {\it Proof.} We know {\rm  (ii)} $\Longleftrightarrow$ {\rm  (iii)} by Lemma \ref {2.3}. {\rm (ii)}$\Longrightarrow$ {\rm (i)} is clear.
{\rm (ii)} $\Longleftarrow$ {\rm (i)}: Assume $\mid u\mid=k$, then $\exists y_{i_{1}},\ldots,y_{i_{k}}\in d(u)$ such that $<y_{i_{1}}\cdots y_{i_{k}},u>\in F^*$ and  $<y_{i_{1}}\cdots y_{i_{k}}, uv>$
$=<y_{i_{1}}\cdots y_{i_{k}}, u>v\neq0$.
 \hfill $\Box$

\begin {Lemma}  \label{p2.3'} If $u$ and $ v$ are two homogeneous elements  with   $\mu (u)\cap \mu (v)=\emptyset$, then  $uv =0$ implies $u=0$ or $v=0$.
\end {Lemma}
\noindent {\it Proof.} Without lost general, there exists $i_0$ such that $\mu (v) \subseteq \{ 1, 2, \cdots i_0\}$ and $\mu (u) \subseteq \{ i_0+1, i_0+2, \cdots n\}$.
Assume $u \not=0$ and $v\not=0$ with $u = \sum \limits_{i=1} ^ s k_i u_i$  and $v = \sum \limits_{j=1} ^ t k_i' v_i$, where $k _i \not=0,$ $k _j' \not=0, $  $u_i  \in B_{1, i_0}$,  $v_j \in B_{i_0+1, n}$ for $1\le i \le i_0, $ $1+i_0\le j  \le n. $ Consequently, $uv = \sum \limits _{i=1} ^ s \sum \limits _{j=1} ^ t k_ik_j' u_iv_j \not= 0$ with $u_iv_j \in B_{1, n}$ and $k_ik_j \not=0$ for $1\le i \le s, $ $1\le j \le t. $  \hfill $\Box$

\section {Relations between    graphes  and Nichols Braided Lie Algebras}

In this section we give relations between connected components of    graphes  and Nichols Braided Lie Algebras.

Obviously,  every pure generalized Dynkin diagram is a graph in terms of  graph theory (see \cite {Ha69}). Conversely,  Assume that  $\Gamma$ is a  graph in terms of  graph theory with
vertex $\{1,  2,  \cdots,  n\}$. We define a matrix  $(q_{ij})_{n \times n}$ as follows:
$q_{ij}q_{ji} \neq 1$ if and only if there exists an edge between $i$ and $j.$ Let $V$ be the braided vector space with braiding matrix $(q_{ij})_{n \times n}$. $\Gamma$ is a
pure generalized Dynkin diagram of $V.$

\begin {Corollary}  \label{5.1}

(i) $\Gamma (u)$ is a connected component of $\Gamma (V)$ with a non-zero monomial $u$ if and only if $\mu (u)$ is a maximal element in $\{\mu (v) \mid v \hbox { is a non-zero monomial in }  \mathfrak L(V) \}$ under order $\subseteq$.

(ii) The following conditions are equivalent.

 (a). $\Gamma (V)$ is connected

 (b). $x_nx_{n-1} \cdots x_1  \in \mathfrak{L}(V)$.

 (c).  $x_1x_{2} \cdots x_n  \in \mathfrak{L}(V)$

(d).  there exists a non-zero monomial $u \in \mathfrak{L}(V)$ with $\mid \mu (u) \mid = \dim V$.

\end {Corollary}
\noindent {\it Proof.}  (i) It follows from Corollary \ref {2.5''}.

(ii) Considering Lemma \ref {p2.3'}
 we have $x_nx_{n-1} \cdots x_1 \not= 0$ and $x_1x_2 \cdots x_n \not= 0$.
Using Corollary \ref {2.5''} we complete the proof.  \hfill $\Box$


\begin{thebibliography}{BD99}



%\bibitem [AS02]{AS02} N. Andruskiewitsch,    H.-J. Schneider,  Pointed Hopf algebras. In New Directions in Hopf
%Algebras,  vol.43 ser MSRI Publications. Cambridge University Press(2002).


%\bibitem[ARS95] {ARS95}M. Auslander,  I. Reiten and S.O. Smal$\phi$,  Representation
%theory of Artin algebras,  Cambridge University Press,  1995.


\bibitem [AS10]{AS10} N. Andruskiewitsch,    H.-J. Schneider,     On the classification of finite-dimensional pointed Hopf algebras,
 Ann. Math.   {\bf 171}  (2010),    375-417.

\bibitem [AHS08] {AHS08} N. Andruskiewitsch,       I. Heckenberger and   H.J. Schneider,
  The Nichols algebra of a semisimple Yetter-Drinfeld module,
 Amer. J. Math. {\bf 132}  (2010),      1493-1547.

%\bibitem[An11] {An11} I. Angiono,        Nichols algebras of unidentified diagonal type,     arXiv:1108.5157.



\bibitem [BFM96]{BFM96} Y. Bahturin,  D. Fishman and S. Montgomery,  On the generalized
Lie structure of associative algebras,  J. Alg. {\bf 96} (1996),  27-48.

\bibitem [BFM01]{BFM01} Y. Bahturin,  D. Fischman and  S. Montgomery,
Bicharacter,  twistings and Scheunert's theorem for Hopf algebra,
J. Alg. {\bf 236} (2001),  246-276.

\bibitem [BMZP92] {BMZP92} Y. Bahturin,  D. Mikhalev,  M. Zaicev and V. Petrogradsky,
Infinite dimensional Lie superalgebras,  Walter de Gruyter Publ.
Berlin,  New York,  1992.


%\bibitem [GM03]{GM03} X. Gomez and S. Majid,   Braided Lie algebras and bicovariant
%differential calculi over coquasitriangular Hopf algebras,
%J. Alg. {\bf 261}(2003),  334--388.

\bibitem [GRR95]{GRR95} D. Gurevich,  A. Radul and V. Rubtsov,
Noncommutative differential geometry related to the Yang-Baxter
equation,  Zap. Nauchn. Sem. S.-Peterburg Otdel. Mat. Inst.
Steklov. (POMI) {\bf 199 } (1992); translation in J. Math. Sci.
{\bf 77 } (1995),  3051--3062.

\bibitem[Gu86] {Gu86} D. I. Gurevich,  The Yang-Baxter equation and the
generalization of formal Lie theory,  Dokl. Akad. Nauk SSSR,  {\bf
288} (1986),  797--801.
\bibitem[He06a]{He06a} I. Heckenberger,   Classification of arithmetic
root systems,   Adv. Math.  {\bf 220} (2009),   59-124.

\bibitem[He06b]{He06b} I. Heckenberger,   The Weyl-Brandt groupoid of a Nichols algebra
of diagonal type,   Invent. Math. {\bf 164} (2006),   175--188.

\bibitem[He05]{He05} I. Heckenberger,       Nichols algebras of diagonal type and arithmetic root systems,     Habilitation,   2005.

\bibitem[Ha69]{Ha69}  Frank Harary,
Graph Theory,  Addison¨CWesley,  USA,  1969

%\bibitem[HS08]{HS08} I. Heckenberger,   H.-J. Schneider,    {Root systems and Weyl groupoids for  Nichols algebras},   {arXiv:0807.0691}.

\bibitem[Hu78]{Hu78}  Humphreys,  James E. Introduction to Lie Algebras and Representation Theory,  Second printing,  revised. Graduate Texts in Mathematics,  9. Springer-Verlag,  New York,  1978.



\bibitem [Ka77]{Ka77} V. G. Kac,  Lie Superalgebras,   Adv. Math. {\bf 26} (1977),
8-96.

\bibitem   [Kh99] {Kh99} V. K. Kharchenko,    A Quantum analog of the
 poincar$\acute{e}$-Birkhoff-Witt theorem,    Algebra and
 Logic,    {\bf 38} (1999),    259-276

% \bibitem [Kh99b]{Kh99b} V. K. Kharchenko,
%An existence condition for multilinear quantum operations,  J. Alg.
%{\bf 217} (1999),  188--228.

%\bibitem[Ma94b] {Ma94b}S. Majid,  Quantum and braided Lie algebras,  J. Geom. Phys.
%{\bf 13} (1994),  307--356.

%\bibitem[Ma95] {Ma95} S. Majid,  Foundations of quantum group,   Cambradge University
%Press,  1995.

\bibitem [Pa98]{Pa98}  B. Pareigis,  On Lie algebras in the category of
Yetter-Drinfeld modules. Appl. Categ. Structures,   {\bf 6}
(1998),  151--175.

\bibitem [Sc79]{Sc79} M. Scheunert,  Generalized  Lie algebras,  J. Math. Phys.
{\bf 20} (1979), 712-720.


%\bibitem [Wo87]{Wo87} Woronowicz,  S.L.,  Compact matrix pseudogroups,  Commun. Math.
%Phys. 111 (1987),  613-665.


\bibitem[WZZ15a] {WZZ15a} W. Wu,    S. Zhang and   Y.-Z. Zhang,     Relationship between Nichols braided Lie
algebras and Nichols algebras,   J. Lie Theory {\bf 25} (2015),      45-63.

\bibitem[WZZ15b] {WZZ15b} W. Wu,    S. Zhang and   Y.-Z. Zhang,     On  Nichols (braided) Lie algebras,   Int. J. Math. {\bf 26} (2015),   1550082. %Also in arXiv:1409.3769.


%\bibitem[WZZ16] {WZZ16} W. Wu,    S. Zhang and   Y.-Z. Zhang,
% Finiteness of Nichols Algebras and Nichols (Braided) Lie Algebras,
%arXiv:1607.07955.


%\bibitem [Zh93]{Zh93} S.C. Zhang,  The Baer radical of generalized matrix rings,
%in Proc. of the Sixth SIAM Conf. on Parallel Processing for
%Scientific Computing,  pp.546--551,  Norfolk,  Virginia,  1993.
%Eds: R.F. Sincovec,  D.E. Keyes,  M.R. Leuze,  L.R. Petzold,  D.A.  Reed.

%\bibitem[ZZH03] {ZZH03} S.C. Zhang,  Y.Z. Zhang and  Y.Y. Han,
%Duality theorems for infinite braided Hopf algebras,  math.QA/0309007.

\bibitem [ZZ03]{ZZ03} S. C. Zhang and  Y.Z. Zhang,
Braided m-Lie Algebras,   Letters in Mathematical Physics,  {\bf 70} (2004),  155-167.


% \bibitem [ZWTZ]{ZWTZ} S. Zhang,   W. Wu,   Zhengtang Tan and   Y.-Z. Zhang,
%Nichols algebras over classical Weyl groups,   Fomin-Kirillov algebras and Lyndon basis,
%     arXiv:1307.8227



\end{thebibliography}
\end {document}